\newif\ifarxiv
\arxivtrue

\ifarxiv
    \documentclass[11pt]{article}
    \usepackage[margin=1in]{geometry}
\else
    \documentclass[smallextended]{svjour3}
    \smartqed
\fi
\usepackage{algorithm}
\usepackage{algpseudocode}
\usepackage{url}

\makeatletter
\newcommand{\LoadPackageUnlessAcm}[2][]{%
    \@ifclassloaded{acmart}{%
        \PackageInfo{preamble}{Skipping package `#2' because acmart class is used}%
    }{%
        \if\relax\detokenize{#1}\relax
            \usepackage{#2}%
        \else
            \usepackage[#1]{#2}%
        \fi
    }%
}
\makeatother

\newif\ifsiamclass
\newif\ifsvjourclass
\makeatletter
\@ifclassloaded{siamart250211}{\siamclasstrue}{\siamclassfalse}
\@ifclassloaded{svjour3}{\svjourclasstrue}{\svjourclassfalse}
\makeatother

\LoadPackageUnlessAcm{caption}
\LoadPackageUnlessAcm{subcaption}
\ifsiamclass\else
    \ifsvjourclass\else
        \LoadPackageUnlessAcm{authblk}
    \fi
    \LoadPackageUnlessAcm{hyperref}
\fi

\usepackage{amstext}
\usepackage{stmaryrd}
\usepackage{algpseudocode}

\usepackage{units}
\usepackage{cancel}
\usepackage{graphicx}

\usepackage{xcolor}
\usepackage{array}
\usepackage{verbatim}
\usepackage{booktabs}
\usepackage{makecell}
\usepackage{calc}
\usepackage{textcomp}
\usepackage{multirow}
\usepackage{tikz,tikzscale}

\usepackage{pgfplots}
\usepackage{pgfplotstable}
\usepackage{pgfkeys}
\pgfplotsset{compat=1.18}

\definecolor{deepgreen}{RGB}{0,100,0}

\pgfplotsset{
    p2/.style={teal, mark=*},
    p3/.style={orange, mark=triangle*},
    p7/.style={blue, mark=square*},
}

\newcommand{\plotwidth}{0.49\columnwidth}
\newcommand{\plotheight}{0.4\columnwidth}

\pgfplotscreateplotcyclelist{paper markers}{
    {teal,    solid, mark=*},
    {orange,   solid, mark=triangle*},
    {blue,     solid, mark=square*},
    {red,      solid, mark=pentagon*},
    {violet,     solid, mark=otimes*},
    {brown,    solid, mark=halfcircle*},
    {black,   solid, mark=diamond*},
}

\colorlet{BarColorOne}{orange!60!white}
\colorlet{BarColorTwo}{teal!80!white}
\colorlet{BarColorThree}{red!80!black}
\colorlet{BarColorFour}{blue!60!white}
\colorlet{BarColorFive}{brown!70!black}

\pgfplotsset{legend swatch/.style={area legend, draw=none}}

\pgfplotsset{
    paperplot/.style={
            width=\plotwidth,
            height=\plotheight,
            cycle list name=paper markers,
            every axis/.append style={font=\small},
            grid=major,
            thick
        }
}

\pgfplotsset{
    longplot/.style={
            width=0.8\columnwidth,
            height=\plotheight,
            cycle list name=paper markers,
            every axis/.append style={font=\small},
            grid=major,
            thick
        }
}

\usepackage[all]{xy}
\usepackage{mathtools}
\usepackage{placeins}
\usetikzlibrary{calc}
\usetikzlibrary{arrows.meta}
\usetikzlibrary{shadings}
\usepgfplotslibrary{fillbetween}

\usepackage[mode=multiuser,status=draft]{fixme}
\fxsetup{innerlayout=inline}
\fxsetup{targetlayout=colorcb}
\fxusetheme{color}
\FXRegisterAuthor{mw}{envmw}{MW}

\definecolor{apply}{rgb}{0.3,.7,0.}
\definecolor{applyface}{rgb}{.9,.95,0.}
\definecolor{invert}{rgb} {0.5,0.,1.}
\colorlet{invertface}{invert!60!red}

\usepackage[colorinlistoftodos,prependcaption,textsize=small]{todonotes}

\usepackage{colortbl}

\usepackage{booktabs}
\usepackage[font=footnotesize,position=b]{subcaption}
\usepackage{adjustbox}

\makeatletter
\def\pgfplots@getautoplotspec into#1{%
        \begingroup
        \let#1=\pgfutil@empty
        \pgfkeysgetvalue{/pgfplots/cycle multi list/@dim}\pgfplots@cycle@dim
        \let\pgfplots@listindex=\pgfplots@numplots
        \pgfkeysgetvalue{/pgfplots/cycle list set}\pgfplots@listindex@set
        \ifx\pgfplots@listindex@set\pgfutil@empty
        \else
            \c@pgf@counta=\pgfplots@listindex
            \c@pgf@countb=\pgfplots@listindex@set
            \advance\c@pgf@countb by -\c@pgf@counta
            \globaldefs=1\relax
            \edef\setshift{%
                \noexpand\pgfkeys{
                    /pgfplots/cycle list shift=\the\c@pgf@countb,
                    /pgfplots/cycle list set=
                }
            }%
            \setshift%
            \globaldefs=0\relax
        \fi
        \pgfkeysgetvalue{/pgfplots/cycle list shift}\pgfplots@listindex@shift
        \ifx\pgfplots@listindex@shift\pgfutil@empty
        \else
            \c@pgf@counta=\pgfplots@listindex\relax
            \advance\c@pgf@counta by\pgfplots@listindex@shift\relax
            \ifnum\c@pgf@counta<0
                \c@pgf@counta=-\c@pgf@counta
            \fi
            \edef\pgfplots@listindex{\the\c@pgf@counta}%
        \fi
        \ifnum\pgfplots@cycle@dim>0
            \c@pgf@counta=\pgfplots@cycle@dim\relax
            \c@pgf@countb=\pgfplots@listindex\relax
            \advance\c@pgf@counta by-1
            \pgfplotsloop{%
                \ifnum\c@pgf@counta<0
                    \pgfplotsloopcontinuefalse
                \else
                    \pgfplotsloopcontinuetrue
                \fi
            }{%
                \pgfkeysgetvalue{/pgfplots/cycle multi list/@N\the\c@pgf@counta}\pgfplots@cycle@N
                \pgfplotsmathmodint{\c@pgf@countb}{\pgfplots@cycle@N}%
                \divide\c@pgf@countb by \pgfplots@cycle@N\relax
                \expandafter\pgfplots@getautoplotspec@
                \csname pgfp@cyclist@/pgfplots/cycle multi list/@list\the\c@pgf@counta @\endcsname
                {\pgfplots@cycle@N}%
                {\pgfmathresult}%
                \t@pgfplots@toka=\expandafter{#1,}%
                \t@pgfplots@tokb=\expandafter{\pgfplotsretval}%
                \edef#1{\the\t@pgfplots@toka\the\t@pgfplots@tokb}%
                \advance\c@pgf@counta by-1
            }%
        \else
            \pgfplotslistsize\autoplotspeclist\to\c@pgf@countd

            \pgfplots@getautoplotspec@{\autoplotspeclist}{\c@pgf@countd}{\pgfplots@listindex}%
            \let#1=\pgfplotsretval
        \fi
        \pgfmath@smuggleone#1%
        \endgroup
    }

\pgfplotsset{
    cycle list set/.initial=
}
\makeatother

\usepackage{etoolbox}
\ifsiamclass\else
\fi
\AtBeginEnvironment{tabular}{\footnotesize}
\AtBeginEnvironment{tabular*}{\footnotesize}
\AtBeginEnvironment{tabularx}{\footnotesize}
\usepackage{lipsum}
\usepackage{amsfonts}
\usepackage{graphicx}
\usepackage{epstopdf}

\makeatletter
\@ifclassloaded{siamart250211}{%
    \usepackage{amssymb}
    \newsiamthm{assumption}{Assumption}
    \newsiamremark{remark}{Remark}
}{%
\@ifclassloaded{svjour3}{%
    \usepackage{amsmath,amssymb}
    \spnewtheorem{assumption}{Assumption}{\bfseries}{\itshape}
}{%
    \usepackage{amsmath,amssymb,amsthm}

    \theoremstyle{definition}
    
    \theoremstyle{remark}
    
    \theoremstyle{plain}
}}
\makeatother

\ifpdf
    \DeclareGraphicsExtensions{.eps,.pdf,.png,.jpg}
\else
    \DeclareGraphicsExtensions{.eps}
\fi

\newcommand{\secref}[1]{Section~\ref{#1}}

\usepackage{amsopn}

\ifarxiv
    \providecommand{\keywordsname}{Keywords}
    \newcommand{\Keywords}[1]{\par\noindent\textbf{\keywordsname:} #1}
    \newcommand{\AMS}[1]{\par\noindent\textbf{AMS subject classifications:} #1}
\else
    \journalname{Journal of Scientific Computing}
\fi

\title{Matrix-Free Local Multigrid for $C^1$ Hermite Discretizations on GPUs}
\begin{document}

\ifarxiv
    \author{Michał Wichrowski}
    \footnotetext{Interdisziplinäres Zentrum für Wissenschaftliches Rechnen (IWR), Ruprecht-Karls-Universität Heidelberg, Germany, \texttt{mwichro@mimuw.edu.pl}}
    \date{}
    \maketitle

    \begin{abstract}

I present a matrix-free local multigrid method for $C^1$ tensor-product Hermite discretizations on adaptively
refined Cartesian meshes, using the clamped biharmonic problem as the model problem. Extending a cell-wise primal--dual formulation
to Hermite elements, I enforce continuity through a constraint operator that is never formed. At refinement edges,
geometric prolongation supplies both the trace and its normal derivative, so the local V-cycle needs no edge operator
and no hanging-node machinery. For the cubic element, an even--odd decomposition based on the signed reflection
symmetry of the Hermite basis, applied along two axes, splits the element matrix into four dense $16\times16$ blocks and reaches $85\%$
of peak memory bandwidth. Chebyshev--Jacobi and vertex-patch smoothers both preserve conformity by construction. In
three dimensions on one NVIDIA A100 in double precision, conjugate gradients preconditioned by the Chebyshev--Jacobi
V-cycle need five or six iterations on uniform meshes and solve at up to $4.23$ million cells per second. On locally
refined meshes the count rises to at most sixteen, ordered by the fraction of cells on a refinement edge and nearly
equal for both smoothers.
     \end{abstract}

    \AMS{65N30, 65N55, 65F08, 65Y20}
\fi

\section{Introduction}
\label{sec:biharmonic_introduction}

The biharmonic equation is a model fourth-order elliptic problem arising in thin-plate bending~\cite{Ciarlet97},
with related formulations in strain-gradient elasticity~\cite{EngelGarikipatiHughesEtAl02}. With clamped boundary
conditions, its natural variational space is $H_0^2(\Omega)$. A conforming finite element approximation must therefore
maintain continuity of the discrete field and its first derivatives across element interfaces. Classical $C^1$
elements meet this requirement through value and derivative degrees of freedom~\cite{Ciarlet78}.

Classical conforming constructions include the triangular plate elements of Argyris, Fried, and
Scharpf~\cite{ArgyrisFriedScharpf68} and the rectangular
element of Bogner, Fox, and Schmit~\cite{BognerFoxSchmit65}. Douglas et
al.~\cite{DouglasDupontPercellEtAl79} developed families of $C^1$ elements with optimal approximation properties
for second- and fourth-order problems.
On Cartesian meshes, tensor-product Hermite interpolation provides a particularly direct construction: the value
trace and the normal-derivative trace are represented by two adjacent layers of element data.

A second line of development relaxes the interelement continuity requirements. Engel et
al.~\cite{EngelGarikipatiHughesEtAl02} studied continuous/discontinuous formulations for fourth-order problems in
structural and continuum mechanics. The $C^0$ interior penalty methods of Brenner and Sung~\cite{BrennerSung05},
reviewed further in~\cite{Brenner11}, retain standard continuous finite element spaces and impose normal-derivative
continuity weakly through consistency and penalty terms on facets. This approach belongs to the wider development of
interior penalty and discontinuous Galerkin methods~\cite{arnold1982interior,arnold2002unified}. The conforming
Hermite formulation considered here uses only volume integrals, with derivative continuity enforced through the
discrete space and preserved by the solver.

Zhang~\cite{Zhang89} is a direct mathematical predecessor for multigrid applied to $C^1$ biharmonic finite elements.
His analysis addresses nonnested spaces, principally the two-dimensional Hsieh--Clough--Tocher macro-element,
using nodal interpolation between uniformly refined levels. With sufficiently many Jacobi-like smoothing steps,
he proves level-independent convergence of W-cycles under the stated regularity and approximation assumptions.
He also explicitly observes that standard multigrid theory applies to the nested bicubic Bogner--Fox--Schmit
rectangular element. Our tensor-product Hermite construction follows this nested setting: geometric prolongation
represents the same coarse polynomial on its children. We address its cell-wise realization, adaptive value and
derivative constraints, and GPU execution.

Zhao~\cite{Zhao05} subsequently analyzed V- and F-cycle methods using the Hsieh--Clough--Tocher element.
For interior penalty discretizations, Brenner and
Zhao~\cite{brenner2005convergence} studied multigrid convergence, while Brenner and Sung~\cite{BrennerSung06}
established convergence estimates for V-, F-, and W-cycles with sufficiently many smoothing steps. Their analysis
uses a smoothing scheme preconditioned by a Poisson solve.

Subspace correction provides a common framework for multigrid and Schwarz iterations~\cite{Xu92}. For the
biharmonic problem, Brenner and Wang~\cite{BrennerWang05} analyzed two-level additive Schwarz preconditioners for
$C^0$ interior penalty methods, and Feng and Karakashian~\cite{FengKarakashian05} studied non-overlapping
two-level Schwarz preconditioners for a discontinuous Galerkin approximation. A complementary connection was
established by Kanschat and Sharma~\cite{KanschatSharma14}: in two dimensions, discrete stream functions relate
divergence-conforming Stokes discretizations to $C^0$ interior penalty biharmonic discretizations. This connection
allows preconditioning ideas to pass between the two problems and helps explain the role of vertex-patch relaxation
in the biharmonic solver of Witte et al.~\cite{witte2025tensor}.

Spline spaces offer another route to $H^2$ conformity through isogeometric analysis. Cho, Pavarino, and
Scacchi~\cite{ChoPavarinoScacchi18} constructed overlapping Schwarz preconditioners for the biharmonic Dirichlet
problem and analyzed their dependence on subdomain and overlap sizes. Sogn and Takacs~\cite{sogn2019robust}
developed multigrid methods with Gauss--Seidel and mass-based smoothing for spline discretizations, proving
robustness with respect to mesh size and, for the mass smoother, spline degree. These studies demonstrate that
strong continuity is compatible with effective multilevel solvers, providing a comparison for Hermite elements
with explicit endpoint value and derivative data.

Matrix-free finite element methods replace stored sparse matrices by evaluations of the bilinear form. Their
efficiency on tensor-product elements rests on sum factorization and fast quadrature, developed in the spectral
and $hp$ finite element literature~\cite{Orszag1980,Melenk99fullydiscrete}. Kronbichler and
Kormann~\cite{Kronbichler2012,Kronbichler2017a} developed cell-based operator evaluation and fast matrix-free
discontinuous Galerkin operators, while Kronbichler and Ljungkvist~\cite{kronbichler2019multigrid} studied complete
matrix-free multigrid solvers on graphics processors. Weber, Kronbichler, and Kreiss~\cite{weber2026hermite}
also developed matrix-free sum-factorized evaluation for conforming Hermite elements, with applications to wave
propagation. Efficient solves require smoothing and inter-grid transfer
to preserve the operator's advantages in data movement, local arithmetic, and parallelism.

Tensor-product local solvers connect Schwarz methods to matrix-free computation. Fast diagonalization originates
with Lynch, Rice, and Thomas~\cite{Lynch1964}; it reduces separable local problems to one-dimensional generalized
eigenproblems and diagonal operations. Overlapping Schwarz methods for spectral elements~\cite{fischer2000overlapping}
and hybrid multigrid/Schwarz algorithms~\cite{LottesFischer05} use local subdomain solves within a multilevel
iteration. More recent developments include fast tensor-product Schwarz smoothers for high-order discontinuous
Galerkin methods~\cite{WitteArndtKanschat21}, approximate tensor-product preconditioners~\cite{pazner2018approximate},
and scalable vertex-star relaxation for high-order finite elements~\cite{brubeck2021scalable}. The local correction
space determines both the couplings resolved in one update and the cost of the local inverse.

On GPUs, Cui et al.~\cite{cui2025implementation} studied the execution of tensor-product patch smoothers, and
Cui and Kanschat~\cite{cui2025multilevel} developed multilevel interior penalty methods.
Witte et al.~\cite{witte2025tensor} address the $C^0$ interior penalty
discretization of the biharmonic problem. Their multigrid method combines overlapping vertex-patch smoothers with
approximate local inverses obtained by fast diagonalization~\cite{Lynch1964}, exploiting tensor structure to make the
patch problems suitable for GPU execution. Their use of an approximate inverse is significant for fourth-order
operators, whose mixed-derivative terms complicate a direct application of a single separable inverse. This is the
closest GPU implementation comparison for the present study, with differences in discretization, continuity
treatment, and local correction spaces.

Local refinement introduces a further requirement: the solver must preserve conformity where cells of different sizes
meet. For Hermite elements, fine-side constraints act on both the field trace and its physical normal derivative,
and the transfer of reference derivative data includes the scale factor of the child map. Local multigrid confines
smoothing to the region present on each refinement level, following the approach developed for adaptive hierarchies
in~\cite{brandt1977mlat,mccormick1986fac,JanssenKanschat11}. Matrix-free geometric multigrid on adaptive meshes
has also been developed in a general parallel framework by Clevenger et al.~\cite{clevenger2021flexible}.
For Hermite elements, both trace layers must be treated consistently during smoothing and coarse-grid correction
while retaining element-local tensor operations.

In this paper, we develop a matrix-free local multigrid method for $C^1$ tensor-product Hermite discretizations on
adaptively refined Cartesian hexahedral meshes, with the clamped biharmonic problem as the model problem. We extend the cell-wise primal--dual formulation of~\cite{wichrowski2026DSS}
to tensor-product Hermite elements. Primal vectors represent conforming fields through element-local copies of their
degrees of freedom, whereas dual vectors retain unassembled residual contributions. Their pairing permits Krylov
iteration without forming a vector of unique global degrees of freedom. Within this formulation, the smoother
converts a dual residual into a conforming primal correction. We describe the two-layer refinement-edge mask and the
geometric transfers that supply the constrained trace data, including the derivative scaling. This separates the
continuity requirement from the cell-wise operator evaluation and expresses adaptive multigrid in terms of local
operations and transfers between levels.

We also study how the algebraic structure of the cubic Hermite element affects GPU execution. Alongside sum
factorization and fast-diagonalization representations, we derive even--odd decompositions based on reflection
symmetry and compare the resulting operator evaluations. We then consider two smoother families: Chebyshev--Jacobi
with synchronization of the value and derivative layers, and a multiplicative vertex-block smoother. The latter
collects contributions from the eight cells incident to a vertex and solves for the eight Hermite degrees of freedom
owned by that vertex. Its correction space is thus smaller than the overlapping Schwarz patch spaces used
in~\cite{witte2025tensor}. This distinction allows us to examine the cost and convergence of a compact local solve
alongside a polynomial smoother within the same conforming discretization and multigrid hierarchy.

Most of the construction depends on the discrete space rather than on the bilinear form. The constraint operator, the
refinement-edge mask, and the transfers are determined by the $C^1$ Hermite space alone. The even--odd decomposition
requires only that the element operator commute with the signed reflection in each coordinate direction. This holds for
every operator assembled from tensor products of the one-dimensional mass, first-derivative, and second-derivative
matrices, including the Hessian form of the plate energy, anisotropic cells, and the mass-plus-biharmonic operator of an
implicit time step. Both smoothers are built from the element operator; only the diagonal, the patch matrix, and the
spectral interval must be recomputed for another form. Kirchhoff plate models~\cite{Ciarlet97} and strain-gradient
elasticity~\cite{EngelGarikipatiHughesEtAl02} lead to the same conformity requirement, as do implicit time steps of
fourth-order evolution equations. The numerical study covers only the biharmonic problem with cubic elements.

The numerical study uses cubic elements in three dimensions and double precision on one NVIDIA A100 GPU. On uniform
meshes, the measured $L^2$ error decreases at fourth order, and conjugate gradients preconditioned by the
Chebyshev--Jacobi V-cycle require five or six iterations over the tested sequence. The two-axis even--odd operator
evaluation reaches an effective memory bandwidth of $85\%$ of the device peak, and the fastest measured uniform-mesh
solve reaches $4.23$ million cells per second. On locally refined meshes, iteration counts depend on the refinement
region, with the largest increase occurring at the first local refinement. These experiments assess discretization
accuracy, solver convergence, and hardware efficiency separately; the adaptive results characterize the tested
refinement sequences rather than establish a mesh-independent convergence bound.

Section~\ref{sec:biharmonic_discretization} introduces the model problem, the Hermite space, and its cell-wise
representation. Section~\ref{sec:biharmonic_multigrid} develops the adaptive local multigrid cycle and its transfer
operators. Section~\ref{sec:biharmonic_operator} derives and compares the matrix-free operator evaluations, and
Section~\ref{sec:biharmonic_smoothing} presents the two smoother families and their costs.
Section~\ref{sec:biharmonic_results} reports discretization and solver results on uniform and locally refined meshes,
followed by conclusions in Section~\ref{sec:biharmonic_conclusion}.

\section{Biharmonic Problem and $H^2$-conforming Discretization}
\label{sec:biharmonic_discretization}

Let $\Omega \subset \mathbb{R}^d$ be a bounded Cartesian domain. We consider the clamped biharmonic problem: find $u
    \in H_0^2(\Omega)$ such that
\begin{equation}
    a(u,v) := \int_\Omega \Delta u\,\Delta v\,dx = \int_\Omega f v\,dx
    \qquad \text{for all } v \in H_0^2(\Omega).
    \label{eq:biharmonic_weak_form}
\end{equation}
Here $H_0^2(\Omega)$ incorporates the boundary conditions $u=0$ and $\partial_nu=0$ on $\partial\Omega$. On
$H_0^2(\Omega)$, the form $a$ coincides with the Hessian inner product $\int_\Omega D^2u:D^2v\,dx$.

We use $d=3$ throughout the algorithms and numerical study, so that $\mathcal{T}_\ell$ is a hierarchy of axis-aligned
hexahedral meshes generated by dyadic refinement. On a mesh without nonmatching interfaces, define
\begin{equation}
    V_\ell := \left\{v \in H_0^2(\Omega): v|_K \circ F_K \in \mathbb{Q}_p([0,1]^d) \text{ for every } K \in \mathcal{T}_\ell\right\},
    \qquad p \geq 3,
    \label{eq:hermite_space}
\end{equation}
where $F_K$ is the affine map from the reference cube to $K$. The $C^1$ trace conditions are part of the definition: on every matching facet, the polynomial trace and all physical first derivatives have identical values from its adjacent cells.

\subsection{Cell-Wise Representation and the Constraint Operator}
\label{sec:biharmonic_cellwise}

Every field is represented cell by cell: $V_{\mathrm{cell},\ell}$ holds one independent block of local degrees of
freedom (DoFs) per cell, and the assembled vector of unique DoFs is never formed. Let $V_{C^1,\ell}\subset
    V_{\mathrm{cell},\ell}$ be the subspace of blocks that represent a globally $C^1$ field on the active mesh, and let the
gather operator $\mathcal{G}_\ell$ replicate the unique DoFs of that field into its cell-local copies, so that
$V_{C^1,\ell}=\operatorname{Im}\mathcal{G}_\ell$.

Following~\cite{wichrowski2026DSS}, continuity is expressed by a constraint operator
\begin{equation}
    C_\ell\colon V_{\mathrm{cell},\ell}\to V_{\mathrm{cell},\ell},
    \qquad C_\ell^2=C_\ell,
    \qquad \operatorname{Im}C_\ell=V_{C^1,\ell},
    \label{eq:biharmonic_constraint_operator}
\end{equation}
so that a cell-wise vector is an admissible primal field exactly when it is a fixed point,
\begin{equation}
    u_\ell=C_\ell u_\ell.
    \label{eq:biharmonic_constraint}
\end{equation}
Idempotence and the range are the only properties used. Which conditions $C_\ell$ encodes is a property of the
discretization, not of the algorithms below: it fixes the clamped boundary data, equates the trace data of adjacent
cells across a matching facet, and slaves the dependent fine-side data at a refinement edge to its coarse-side support.
The fourth-order setting differs from a nodal $C^0$ one in that the constrained trace datum is the pair
$(u,\partial_nu)$ rather than the value alone. We never form $C_\ell$, and no explicit hanging-node machinery appears
in the solver.

Dual vectors $\tilde r_\ell\in V_{\mathrm{cell},\ell}^*$ are raw unassembled element contributions, meaningful without
any synchronization; their assembled representation is $C_\ell^T\tilde r_\ell$. For a primal field
satisfying~\eqref{eq:biharmonic_constraint}, the Euclidean pairing is therefore already the global one,
\begin{equation}
    \langle\tilde r_\ell,u_\ell\rangle
    =\langle\tilde r_\ell,C_\ell u_\ell\rangle
    =\langle C_\ell^T\tilde r_\ell,u_\ell\rangle,
    \label{eq:biharmonic_pairing}
\end{equation}
and the discrete problem reads: find $u_\ell=C_\ell u_\ell$ with
\begin{equation}
    C_\ell^TA_\ell C_\ell u_\ell=C_\ell^T\tilde b_\ell,
    \label{eq:biharmonic_constrained_system}
\end{equation}
which is equivalent to the assembled system on the $C^1$ space.

Operator application, vector updates, and inter-grid transfers act cell by cell on this storage. The one component that
converts a dual residual into a primal correction is the preconditioner
\begin{equation}
    P_\ell\colon V_{\mathrm{cell},\ell}^*\to V_{C^1,\ell},
    \label{eq:biharmonic_preconditioner_requirement}
\end{equation}
and with it the multigrid smoother. Its output must be a fixed point of $C_\ell$; how it achieves this is the subject
of \secref{sec:biharmonic_smoothing}. Nothing else in the solver has to know how continuity is defined.

Let $\mathcal{P}_{\ell-1}^\ell$ be the geometric prolongation given by the embedding $V_{\ell-1}\subset V_\ell$. Beyond
\eqref{eq:biharmonic_constraint_operator}, the multigrid construction uses one property of the constraint operator,
\begin{equation}
    C_{\ell}\mathcal{P}_{\ell-1}^\ell=\mathcal{P}_{\ell-1}^\ell.
    \label{eq:biharmonic_constraint_prolongation}
\end{equation}
A prolongated coarse field is a $C^1$ field: the embedding already places the interpolated trace data on a refinement
edge, leaving nothing for the constraints to enforce. Any formulation whose constraint operator
satisfies~\eqref{eq:biharmonic_constraint_prolongation} is covered by what follows.

\subsection{Tensor-Product Hermite Element}
\label{sec:biharmonic_hermite}

The element used in this work makes the constrained data of $C_\ell$ explicit. Let $\{\psi_i\}_{i=0}^p \subset
    \mathbb{P}_p([0,1])$ be a one-dimensional Hermite basis whose endpoint functionals satisfy
\begin{equation}
    q(0)=q_0, \qquad q'(0)=q_1, \qquad q'(1)=q_{p-1}, \qquad q(1)=q_p,
    \label{eq:biharmonic_hermite_endpoints}
\end{equation}
for $q=\sum_{i=0}^p q_i\psi_i$. The remaining basis functions have vanishing value and derivative at both endpoints.
The endpoint derivative DoFs are reference derivatives; on a cell of width $h_\alpha$ their associated physical
derivatives are $h_\alpha^{-1}q_1$ and $h_\alpha^{-1}q_{p-1}$.

The tensor-product basis $\psi_{i_1}(\hat x_1)\psi_{i_2}(\hat x_2)\psi_{i_3}(\hat x_3)$ separates the DoFs that
determine a facet trace from those that determine its normal derivative. On the reference facet $\hat x_1=0$, for
example, the DoF layer $i_1=0$ determines $u|_F$, whereas the layer $i_1=1$ determines $\partial_{\hat x_1}u|_F$.
Tangential derivatives are obtained by differentiating the trace polynomial. Equality of these two layers, after
physical scaling, is therefore sufficient for equality of the field and its full gradient on a matching facet, and the
conditions collected in $C_\ell$ act on the same two layers: they equate them across a matching facet, prescribe them
by interpolation from the coarse side at a refinement edge, and fix them on a clamped boundary facet. At edges and
vertices shared by several boundary facets, the tensor-product representation applies all compatible conditions to the
same local DoFs.

The prolongation $\mathcal{P}_{\ell-1}^\ell$ is the tensor product of one-dimensional child interpolation matrices,
whose derivative rows carry the chain-rule factor of the child map: a reference derivative DoF is scaled by the ratio
of child to parent cell width when it is transferred. This scaling is the one point at which the transfer differs from
a value-only $C^0$ construction.

\section{Local Multigrid}
\label{sec:biharmonic_multigrid}

\subsection{Cell-Wise Inter-Grid Transfers}
\label{sec:biharmonic_transfer_operators}

For a parent cell $K$ and its child $K_c$, define the local prolongation by
\begin{equation}
    P_{K\to K_c}=P_{c_1}\otimes P_{c_2}\otimes P_{c_3},
    \qquad
    \mathcal{R}_{\mathrm{cell},\ell}^{\ell-1}=
    \bigl(\mathcal{P}_{\mathrm{cell},\ell-1}^{\ell}\bigr)^T,
    \label{eq:biharmonic_transfers}
\end{equation}
where each $P_{c_\alpha}$ evaluates a parent Hermite polynomial on child $c_\alpha$, with the derivative scaling of
\secref{sec:biharmonic_hermite}. Prolongation acts on primal DoFs and restriction accumulates raw dual contributions
from all children.

The constraint identity~\eqref{eq:biharmonic_constraint_prolongation} yields the two relations required by local
multigrid:
\begin{align}
    \mathcal{P}_{\mathrm{cell},\ell-1}^{\ell}\mathcal{G}_{\ell-1}
     & =\mathcal{G}_{\ell}\mathcal{P}_{\ell-1}^{\ell},
    \label{eq:biharmonic_prolongation_commutes}                    \\
    \mathcal{G}_{\ell-1}^T\mathcal{R}_{\mathrm{cell},\ell}^{\ell-1}
     & =\mathcal{R}_{\ell}^{\ell-1}C_{\ell}^T\mathcal{G}_{\ell}^T.
    \label{eq:biharmonic_restriction_commutes}
\end{align}
The first relation says that cell-wise prolongation yields a $C^1$ field. The second says that restricting an
unassembled residual implicitly applies $C_\ell^T$, so no constraint matrix is needed on the dual side either. Both
follow from the duality argument of~\cite{wichrowski2026DSS}, which uses only the fixed-point
property~\eqref{eq:biharmonic_constraint} and the prolongation
identity~\eqref{eq:biharmonic_constraint_prolongation}.

\subsection{Level-Wise Decomposition of an Adaptive Mesh}
\label{sec:biharmonic_level_decomposition}

Adaptive refinement subdivides a cell into $2^d$ children and may be applied recursively, so the discrete problem is
posed on the leaves of the refinement tree, the active mesh $\mathcal{T}_{\mathrm{act}}$. Adjacent leaves may belong to
different levels, and on the fine side of a refinement edge the value and normal-derivative data are constrained by the
coarse side through $C_\ell$. A classical V-cycle would require a sequence of uniformly refined meshes, which is
unavailable here: embedding $\mathcal{T}_{\mathrm{act}}$ into a notional uniform fine mesh would introduce DoFs the
discretization does not carry and would make the per-level cost grow with the uniform mesh instead of with the refined
region. Local multigrid~\cite{brandt1977mlat,mccormick1986fac,JanssenKanschat11} avoids this by smoothing on level
$\ell$ only within the region covered by cells of that level.

Let $\mathcal{T}_\ell$ be the collection of all cells of level $\ell$, active or further refined, and let
$V_{\mathrm{cell},\ell}$ be the associated cell-wise space. Its active subspace $V_{\mathrm{cell},\ell}^{\mathrm{act}}$
consists of the cell-wise fields vanishing on every refined cell of $\mathcal{T}_\ell$. Every active cell has a unique
level, so
\begin{equation}
    V_{\mathrm{cell},\mathrm{act}}=\bigoplus_{\ell=0}^{L}V_{\mathrm{cell},\ell}^{\mathrm{act}}.
    \label{eq:biharmonic_active_decomposition}
\end{equation}
This splitting is the interface between the active mesh and the multigrid hierarchy. The unassembled active residual
$\tilde r_{\mathrm{act}}$ is distributed by restricting it to each summand and extending it by zero over the refined
cells of that level, giving the level right-hand sides $\tilde r_\ell\in V_{\mathrm{cell},\ell}^*$; conversely, the
level corrections $e_\ell$ are collected by summing their components on active cells. The work of one cycle is
therefore proportional to the number of active cells summed over the levels.

\subsection{Refinement-Edge Masking}
\label{sec:biharmonic_edge_masking}

Let $E_\ell$ contain every level-$\ell$ DoF whose value or first derivative is prescribed by a coarser active neighbor,
together with the DoFs fixed by the clamped boundary conditions, and let $S_\ell$ be its complement, the interior of
the region covered by the cells of level $\ell$. The smoother updates $S_\ell$ only.

This partition splits the assembled level operator into the blocks $A^{SS}$, $A^{SE}$, $A^{ES}$, and $A^{EE}$. In the classical
formulation the smoother solves $A^{SS}u^S = f^S-A^{SE}u^E$ with the prolongated edge data $u^E$ as Dirichlet values,
and a separate edge operator evaluates $r^E=f^E-A^{ES}u^S-A^{EE}u^E$ to drive the coarse-grid correction. Neither block
is formed here. The full volumetric operator is evaluated on every cell of $\mathcal{T}_\ell$ in one pass; it reads the
edge data $u^E$ during the cell integration and thereby applies $A^{SE}$ on the fly. The smoother masks the dual
residual on $E_\ell$ and extends its primal correction by zero there, which reproduces the interior relaxation, and the
unmasked residual evaluated after smoothing carries the outward fluxes $A^{ES}$ that the plain
restriction~\eqref{eq:biharmonic_transfers} transports to the coarser level. In the $C^1$ setting the mask covers both
trace layers: freezing only the value layer would leave the normal-derivative data of the edge free and destroy
conformity.

\subsection{The Local V-Cycle}
\label{sec:biharmonic_local_vcycle}

Algorithm~\ref{alg:biharmonic_local_vcycle} states one cycle. Its level right-hand side enters already carrying the
component of $\tilde r_{\mathrm{act}}$ placed on that level by~\eqref{eq:biharmonic_active_decomposition}, so the
restricted fine residual is added to it rather than overwriting it. The edge set $E_\ell$ is fixed during both
smoothing stages, and the coarse-grid correction supplies the value and normal-derivative data on $E_\ell$ through
prolongation, so~\eqref{eq:biharmonic_constraint} is maintained without an assembled edge operator or an explicit
constraint matrix.

\begin{algorithm}[H]
    \caption{Local $C^1$ multigrid V-cycle: $\mathrm{MG}_\ell(\tilde b_\ell,u_\ell)$}
    \label{alg:biharmonic_local_vcycle}
    \begin{algorithmic}[1]
        \Require Level $\ell$; cell-wise iterate $u_\ell$; unassembled level right-hand side $\tilde b_\ell$ carrying
        the level-$\ell$ component of $\tilde r_{\mathrm{act}}$.
        \Ensure Updated $u_\ell$.
        \If{$\ell=0$}
        \State \Return coarse-level solve for $\tilde b_0$
        \EndIf
        \State $u_\ell\gets\mathcal{S}_\ell^{\nu_1}(\tilde b_\ell,u_\ell)$ \Comment{Pre-smoothing on $S_\ell$}
        \State $\tilde r_\ell\gets\tilde b_\ell-A_\ell u_\ell$ \Comment{Unassembled cell-wise residual}
        \State $\tilde b_{\ell-1}\gets\tilde b_{\ell-1}+\mathcal{R}_{\mathrm{cell},\ell}^{\ell-1}\tilde r_\ell$
        \Comment{Accumulate, no constraint transpose}
        \State $u_{\ell-1}\gets\mathrm{MG}_{\ell-1}(\tilde b_{\ell-1},u_{\ell-1})$
        \State $u_\ell\gets u_\ell+\mathcal{P}_{\mathrm{cell},\ell-1}^{\ell}u_{\ell-1}$ \Comment{Prolongate; sets the
            $C^1$ data on $E_\ell$}
        \State $u_\ell\gets\mathcal{S}_\ell^{\nu_2}(\tilde b_\ell,u_\ell)$ \Comment{Post-smoothing on $S_\ell$}
        \State \Return $u_\ell$
    \end{algorithmic}
\end{algorithm}

Relations~\eqref{eq:biharmonic_prolongation_commutes} and~\eqref{eq:biharmonic_restriction_commutes} identify the
transfers of Algorithm~\ref{alg:biharmonic_local_vcycle} with those of the classical local
V-cycle~\cite{JanssenKanschat11} built from the edge constraints in $C_\ell$. If the smoother of a level maps a $C^1$
iterate to the gathered image of a classical local smoother on $S_\ell$, then the whole cycle agrees with the classical
one after gathering. Both smoothers of \secref{sec:biharmonic_smoothing} have this property: each is a classical
smoother of the assembled level operator, returns a $C^1$ correction, and leaves $E_\ell$ untouched. The equivalence is
algebraic; we study the convergence of the cycle numerically in \secref{sec:biharmonic_results}.

\section{Evaluation of the Biharmonic Operator}
\label{sec:biharmonic_operator}

This section describes alternative matrix-free evaluations of the same three-dimensional biharmonic operator. The
methods differ only in the realization of the local tensor contractions; they use the same Hermite DoFs, boundary
treatment, and cell traversal.

\subsection{Element Operator and Sum Factorization}
\label{sec:biharmonic_sum_factorization}

Let $M$, $K$, and $D$ denote the one-dimensional matrices with entries
\begin{equation}
    M_{ij}=\int_0^1\psi_i\psi_j\,d\hat x, \qquad
    K_{ij}=\int_0^1\psi_i'\psi_j'\,d\hat x, \qquad
    D_{ij}=\int_0^1\psi_i''\psi_j''\,d\hat x.
    \label{eq:biharmonic_1d_factors}
\end{equation}
For an isotropic hexahedral cell of width $h$, the cell matrix of the Hessian form $\int_K D^2u:D^2v\,dx$ is
\begin{equation}
    \begin{aligned}
        A_K = h^{-1}\bigl( & D\otimes M\otimes M + M\otimes D\otimes M + M\otimes M\otimes D             \\
                           & + 2K\otimes K\otimes M + 2K\otimes M\otimes K + 2M\otimes K\otimes K\bigr).
    \end{aligned}
    \label{eq:biharmonic_3d_operator}
\end{equation}
The pure second derivatives give the $D$ terms and the mixed ones $\partial_i\partial_j$, $i\neq j$, give the $K\otimes K$
terms. The cell matrix of $\int_K\Delta u\,\Delta v\,dx$ differs from $A_K$, but the two forms agree on $H_0^2(\Omega)$
(\secref{sec:biharmonic_discretization}), so the assembled operator on $V_\ell$ is that of~\eqref{eq:biharmonic_weak_form}. An anisotropic mesh requires one physical scale factor per derivative direction rather than the isotropic factor shown here.

Sum factorization evaluates~\eqref{eq:biharmonic_3d_operator} through successive one-dimensional contractions, as in
the matrix-free biharmonic multigrid implementation of Witte et al.~\cite{witte2025tensor}. Applied directly to the six
terms, it needs twelve contractions per element vector even when intermediate products are shared, and the terms must
be accumulated before the last contraction. Fast diagonalization, \secref{sec:biharmonic_fast_diagonalization}, keeps
that data flow and the same contraction shapes, but regroups the six terms into two parts that share their first and
last contractions, so that one evaluation issues ten contractions and the two parts are summed inside the last one. We
therefore take fast diagonalization, rather than sum factorization, as the baseline against which the decompositions
below are measured.

\subsection{Fast Diagonalization}
\label{sec:biharmonic_fast_diagonalization}

Let $\mathcal{M}=M\otimes M\otimes M$ and let $\mathcal{L}=K\otimes M\otimes M+M\otimes K\otimes M+M\otimes M\otimes K$
be the element Laplacian. Write the element operator as
\begin{equation}
    \begin{gathered}
        A_K = h^{-1}\bigl(\mathcal{L}\mathcal{M}^{-1}\mathcal{L} + \mathcal{C}\bigr),
        \qquad
        \widetilde D = D - KM^{-1}K, \\
        \mathcal{C} = \widetilde D\otimes M\otimes M + M\otimes\widetilde D\otimes M + M\otimes M\otimes\widetilde D.
    \end{gathered}
    \label{eq:biharmonic_fdm_operator}
\end{equation}
Expanding $\mathcal{L}\mathcal{M}^{-1}\mathcal{L}$, the products of two different terms of $\mathcal{L}$ give exactly
$2(K\otimes K\otimes M+K\otimes M\otimes K+M\otimes K\otimes K)$, whereas the product of a term with itself gives
$KM^{-1}K$ in place of $D$ in one direction. The correction $\mathcal{C}$ replaces these diagonal products by
$D\otimes M\otimes M$ and its permutations, so the split reproduces~\eqref{eq:biharmonic_3d_operator} exactly.

Let $V_L$ and $V_C$ collect the generalized eigenvectors of the pairs $(K,M)$ and $(\widetilde D,M)$, normalized so
that $V^TMV=I$, with eigenvalues $\{\lambda_i^L\}_{i=0}^p$ and $\{\lambda_i^C\}_{i=0}^p$. With the transforms
$Q_L=V_L^TM$ and $Q_C=V_C^TM$, one has $M=Q^TQ$ for both, $K=Q_L^T\operatorname{diag}(\lambda^L)Q_L$, and
$\widetilde D=Q_C^T\operatorname{diag}(\lambda^C)Q_C$. The diagonal tensors are
\begin{equation}
    (\Lambda_L)_{ijk}=\bigl(\lambda_i^L+\lambda_j^L+\lambda_k^L\bigr)^2,
    \qquad
    (\Lambda_C)_{ijk}=\lambda_i^C+\lambda_j^C+\lambda_k^C.
    \label{eq:biharmonic_fdm_diagonals}
\end{equation}
The square in $\Lambda_L$ arises because $\mathcal{M}^{-1}$ cancels one pair of transforms between the two factors of
$\mathcal{L}$. This is the generalized-eigenvalue construction of fast diagonalization~\cite{Lynch1964}, applied once
to each of the two parts.

Their operator evaluations are
\begin{align}
    \mathcal{L}\mathcal{M}^{-1}\mathcal{L}
     & = (Q_L^T\otimes Q_L^T\otimes Q_L^T)\,\Lambda_L\,(Q_L\otimes Q_L\otimes Q_L),
    \label{eq:biharmonic_fdm_l}                                                     \\
    \mathcal{C}
     & = (Q_C^T\otimes Q_C^T\otimes Q_C^T)\,\Lambda_C\,(Q_C\otimes Q_C\otimes Q_C).
    \label{eq:biharmonic_fdm_c}
\end{align}
Thus FDM evaluates each part by three contractions with $Q_L$ or $Q_C$, a diagonal scaling, and three contractions
with the corresponding transpose. It evaluates the operator, not its inverse.

For cubic Hermite elements, $p=3$ and each one-dimensional transform has order four. The two FDM parts have the same
input vectors and their results are added, so their first and final contractions can be stacked without discarding
entries. In the row-oriented realization of~\eqref{eq:biharmonic_fdm_l}--\eqref{eq:biharmonic_fdm_c}, the forward
transforms form $[Q_L^T\mid Q_C^T]\in\mathbb{R}^{4\times8}$ and the backward transforms form
$[Q_L;Q_C]\in\mathbb{R}^{8\times4}$; the latter product also accumulates the two parts. One evaluation therefore
issues ten contractions: the stacked forward one, two forward and two backward contractions in each part, and the
stacked backward one. Applied to an eight-row tile, the forward product has the $8\times8\times4$ tensor-core MMA
shape.

\subsection{Even--Odd Decompositions}
\label{sec:biharmonic_even_odd}

A matrix that commutes with the exchange matrix $J$, which has ones on the antidiagonal, is called centrosymmetric.
Its even--odd block diagonalization is classical~\cite{cantoni1976centrosymmetric}, and it halves the cost of a
product with a symmetric centrosymmetric matrix once the two blocks are stored~\cite{melman2000centrosymmetric}. In
spectral and matrix-free finite element methods, the same splitting is applied to the one-dimensional interpolation
and derivative kernels of sum factorization~\cite{solomonoff1992fast,kopriva2009implementing,Kronbichler2017a}. We
apply it to the dense element matrix of the Hermite element instead. Two points differ: the reflection carries a sign on
the derivative DoFs, and it is applied independently in each coordinate direction.

Order the one-dimensional cubic Hermite DoFs as left value, left derivative, right derivative, and right value.
Reflection interchanges the endpoint values and interchanges the endpoint derivatives with a sign change,
\begin{equation}
    \mathcal{R}(u_0,u_1,u_2,u_3)^T=(u_3,-u_2,-u_1,u_0)^T.
    \label{eq:biharmonic_hermite_reflection}
\end{equation}

The elementary transformation used below replaces a pair of entries by their sum and difference,
\begin{equation}
    (a,b)\longmapsto(a+b,a-b).
    \label{eq:biharmonic_scalar_butterfly}
\end{equation}
The sum is unchanged when $a$ and $b$ are interchanged, whereas the difference changes sign, so the transformation
separates the even and odd parts of the pair. It requires one addition and one subtraction. When its dependencies are
drawn as a data-flow graph, both inputs connect to both outputs and the crossing connections form the characteristic
two-wing shape; for this reason, we refer to the sum--difference transformation as a \emph{butterfly}.

For the endpoint values the mirrored pair is $(u_0,u_3)$. For the endpoint derivatives it is $(u_1,-u_2)$ because
reflection reverses the direction of differentiation. Applying the sum--difference operation to these two pairs gives
\begin{equation}
    u_0^+=u_0+u_3,\qquad u_0^-=u_0-u_3,
    \qquad
    u_1^+=u_1-u_2,\qquad u_1^-=u_1+u_2.
    \label{eq:biharmonic_hermite_butterfly_entries}
\end{equation}
Collecting the even entries first and the odd entries second yields the unnormalized transform
\begin{equation}
    \begin{bmatrix}
        u_0^+ \\ u_1^+\\ u_0^-\\ u_1^-
    \end{bmatrix}
    =H\boldsymbol{u}
    :=
    \begin{bmatrix}
        1 & 0 & 0  & 1  \\
        0 & 1 & -1 & 0  \\
        1 & 0 & 0  & -1 \\
        0 & 1 & 1  & 0
    \end{bmatrix}
    \begin{bmatrix}u_0\\u_1\\u_2\\u_3\end{bmatrix},
    \qquad H^{-1}=\frac{1}{2}H^T.
    \label{eq:biharmonic_hermite_butterfly}
\end{equation}

Thus $(u_0^+,u_1^+)$ are invariant under~\eqref{eq:biharmonic_hermite_reflection}, whereas $(u_0^-,u_1^-)$ change sign.
The inverse butterfly reconstructs each original entry by averaging the sum and difference, which accounts for the
factor $1/2$ in $H^{-1}$. The derivative signs in $H$ are essential: treating the derivative pair like the value pair
would describe index reversal rather than reflection of the Hermite function. Equivalently, $\mathcal{R}=SJS$ with
$S=\operatorname{diag}(1,1,-1,1)$: once the right-endpoint derivative is negated, so that both derivative DoFs point
into the cell, the one-dimensional matrices $M$, $K$, and $D$ are symmetric centrosymmetric.

Let $S\subseteq\{x,y,z\}$ contain the transformed axes, let $r=|S|$, and define
\begin{equation}
    H_S=T_x\otimes T_y\otimes T_z,
    \qquad
    T_\alpha=
    \begin{cases}
        H,   & \alpha\in S,    \\
        I_4, & \alpha\notin S.
    \end{cases}
    \qquad
    H_S^{-1}=2^{-r}H_S^T.
    \label{eq:biharmonic_partial_parity_transform}
\end{equation}
Let $P_S$ permute the transformed vector entries into sector-major order, with the parity tuple first and the retained
local indices second, and set $E_S=P_SH_S$. Its inverse is $E_S^{-1}=2^{-r}H_S^TP_S^T$.

The tensor-product notation does not imply that the dense $64\times64$ matrix $H_S$ is formed or applied. An
$x$-butterfly applies~\eqref{eq:biharmonic_scalar_butterfly} independently to every pair of entries mirrored in $x$,
with the $y$ and $z$ indices held fixed. A subsequent $y$-butterfly applies the same operation to the resulting data
with $x$ parity and $z$ held fixed, and similarly for $z$. Each transformed axis therefore doubles the number of parity
sectors and halves the number of entries in each sector. The inverse evaluates the same stages in reverse order and
contributes one factor $1/2$ per transformed axis.

Since the Cartesian element operator commutes with the signed reflection in every coordinate direction, this partial
transform is an exact block diagonalization,
\begin{equation}
    E_S A_K E_S^{-1}
    =\bigoplus_{\boldsymbol{\sigma}\in\{+,-\}^{r}}A_{\boldsymbol{\sigma}},
    \qquad
    A_{\boldsymbol{\sigma}}\in\mathbb{R}^{n_r\times n_r},
    \qquad
    n_r=2^r4^{3-r}=\frac{64}{2^r}.
    \label{eq:biharmonic_partial_parity_blocks}
\end{equation}
Each transformed axis contributes one parity bit and retains two entries within each parity. Every untransformed
axis retains all four entries. The operator evaluation consists of a forward butterfly, independent dense
products with the diagonal blocks, and the inverse butterfly,
\begin{equation}
    A_K\boldsymbol{u}
    =E_S^{-1}
    \left(\bigoplus_{\boldsymbol{\sigma}\in\{+,-\}^{r}}A_{\boldsymbol{\sigma}}\right)
    E_S\boldsymbol{u}.
    \label{eq:biharmonic_partial_parity_evaluation}
\end{equation}

For one element vector, the useful arithmetic in the dense blocks is
\begin{equation}
    2\,2^r n_r^2
    =2\,2^r\left(\frac{64}{2^r}\right)^2
    =\frac{8192}{2^r}.
    \label{eq:biharmonic_even_odd_flops}
\end{equation}
floating-point operations, counting a fused multiply--add as two operations. For $r=1$ this is the halving
of~\cite{melman2000centrosymmetric}; each further axis halves the count again. This count excludes the butterflies:
each transformed axis adds $64$ additions in the forward transform and $64$ in the inverse transform. The factor
$2^{-r}$ may be folded into the sector matrices.

\begin{table}[htbp]
    \centering
    \caption{Partial even--odd decompositions for a three-dimensional cubic Hermite element. The FLOP count excludes
        butterfly additions.}
    \label{tab:biharmonic_even_odd_variants}
    \begin{tabular}{cccc r}
        \toprule
        Variant & split axes    & sectors & block order & FLOPs \\
        \midrule
        0       & $\varnothing$ & 1       & 64          & 8192  \\
        I       & $x$           & 2       & 32          & 4096  \\
        II      & $x,y$         & 4       & 16          & 2048  \\
        \bottomrule
    \end{tabular}
\end{table}

\paragraph{Variant 0: dense matrix--matrix multiplication.}
No butterfly is applied. For a batch $U\in\mathbb{R}^{64\times b}$ of $b$ element vectors, the evaluation is the single
product $A_KU$. This form has no parity-sector layout, permutations, or butterfly dependencies. It retains every
coupling in the assembled $64\times64$ matrix and therefore performs $8192$ floating-point operations per element
vector.

\paragraph{Variant I: one-axis decomposition.}
Take $S=\{x\}$. The forward $x$-butterfly forms two vectors $\widetilde{\boldsymbol{u}}_+$ and
$\widetilde{\boldsymbol{u}}_-$, each indexed by two internal $x$ modes and all four modes in both $y$ and $z$. Hence
each vector has $2\cdot4\cdot4=32$ entries. The two products $A_+\widetilde{\boldsymbol{u}}_+$ and
$A_-\widetilde{\boldsymbol{u}}_-$ use independent $32\times32$ blocks; the inverse $x$-butterfly then recombines the
two results with the factor $1/2$. This variant performs $4096$ dense-product operations. It has only one forward and
one inverse butterfly level and uses the widest blocks among the decomposed forms, but retains all couplings associated
with reflections in $y$ and $z$.

\paragraph{Variant II: two-axis decomposition.}
Take $S=\{x,y\}$ and leave $z$ unsplit. Successive signed butterflies first separate $x$ parity and then $y$ parity,
producing the four sectors $\widetilde{\boldsymbol{u}}_{++}$, $\widetilde{\boldsymbol{u}}_{+-}$,
$\widetilde{\boldsymbol{u}}_{-+}$, and $\widetilde{\boldsymbol{u}}_{--}$. Within each sector, the two internal modes
from each split axis and the four original $z$ modes give $2\cdot2\cdot4=16$ entries. Four independent $16\times16$
products evaluate the sector operators. The inverse $y$- and $x$-butterflies restore the original entry order and
supply the normalization factor $1/4$.

Each of the four sector products is a dense $16\times16$ product without padding, and the decomposition stops here.
On isotropic cells any pair of axes is equivalent; choosing $x$ and $y$ fixes the pair to the local DoF layout and
leaves the remaining $z$ index contiguous within each sector.

\subsection{Comparison of Operator Evaluations}
\label{sec:biharmonic_operator_comparison}

To compare the efficiency of the variants, we measure their wall time, effective memory bandwidth, and throughput on
the same set of element vectors. The benchmark contains $4\,194\,304$ cells, corresponding to the total of
$268\,435\,456$ vector entries. The layout interleaves batches of $32$ cells in its fastest-running index. All
measurements reported in this work, here and in \secref{sec:biharmonic_smoother_performance} and
\secref{sec:biharmonic_results}, use one NVIDIA A100-SXM4-80GB GPU, whose peaks are $2039$\,GB/s of memory bandwidth,
$9.7$\,TFLOP/s of vector FP64, and $19.5$\,TFLOP/s of tensor-core FP64.

Each evaluation reads and writes $64$ vector entries per cell. We therefore compute the effective bandwidth from the
minimum data volume
\begin{equation}
    V_{\mathrm{mem}}
    =2\cdot64\cdot 4\,194\,304\cdot 8\ \text{bytes}
    =4.295\ \text{GB}.
    \label{eq:biharmonic_operator_benchmark_traffic}
\end{equation}
The small cell-operator tables are reused across all cells and remain outside this streaming byte count. Consequently,
the reported bandwidth measures how rapidly each realization processes the common input and output data, rather than
the total number of arithmetic instructions it issues.

Before timing, the FDM reconstruction is compared with the assembled $64\times64$ cell matrix. The FDM, dense,
one-axis, and two-axis evaluations are then applied to the same random element vectors and compared with the
assembled result. All four evaluations agree with the reference to FP64 roundoff.

\begin{table}[htbp]
    \centering
    \caption{Evaluation of the cubic three-dimensional biharmonic operator.
        Throughput is reported in cells per second. Effective bandwidth uses the minimum read--write volume
        in~\eqref{eq:biharmonic_operator_benchmark_traffic}.}
    \label{tab:biharmonic_operator_performance}
    \begin{tabular}{lrrr}
        \toprule
        Evaluation                  & time [ms]      & bandwidth [GB/s] & throughput [$10^9$ cells/s] \\
        \midrule
        FDM                         & 4.857          & 884.3            & 0.864                       \\
        Dense matrix--matrix        & 3.738          & 1149.1           & 1.122                       \\
        \textbf{Two-axis even--odd} & \textbf{2.475} & \textbf{1735.4}  & \textbf{1.695}              \\
        One-axis even--odd          & 2.532          & 1696.5           & 1.657                       \\
        \bottomrule
    \end{tabular}
\end{table}
\FloatBarrier

The two-axis decomposition gives the highest measured throughput, $1.695\times10^9$ cells/s, and reaches
$1735.4$\,GB/s, or $85\%$ of peak memory bandwidth. It is $1.96\times$ faster than FDM, $1.51\times$ faster than the
dense matrix--matrix evaluation. The one-axis variant is within $2.3\%$ of its throughput although it performs twice
the dense arithmetic: both even--odd variants are limited by memory bandwidth.

\subsection{Roofline of the Operator Realizations}
\label{sec:biharmonic_roofline}

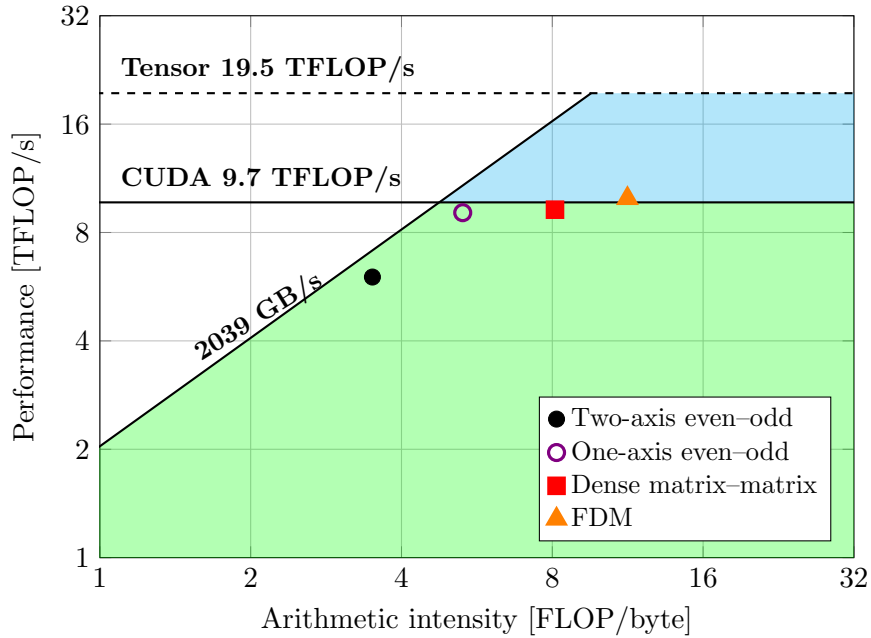
\begin{figure}[htbp]
    \centering
\begin{tikzpicture}
    \begin{loglogaxis}[
            width=0.7\linewidth,
            height=0.53\linewidth,
            xlabel={Arithmetic intensity [FLOP/byte]},
            ylabel={Performance [TFLOP/s]},
            log basis x=2,
            xmin=1, xmax=32,
            xtick={1,2,4,8,16,32},
            xticklabels={1,2,4,8,16,32},
            log basis y=2,
            ymin=1, ymax=32,
            ytick={1,2,4,8,16,32},
            yticklabels={1,2,4,8,16,32},
            grid=both,
            legend pos=south east,
            legend columns=1,
            legend cell align=left,
            legend style={font=\small},
        ]
        \addplot[green, fill=green, draw=none, fill opacity=0.3, forget plot]
        coordinates {(0.001,0.002039) (4.76,9.7) (1000,9.7) (1000,0.0001) (0.001,0.0001)};
        \addplot[cyan, fill=cyan, draw=none, fill opacity=0.3, forget plot]
        coordinates {(4.76,9.7) (9.56,19.5) (1000,19.5) (1000,9.7)};

        \addplot[thick, domain=0.001:9.56, samples=2, forget plot] {2.039*x};
        \addplot[thick, domain=0.001:1000, samples=2, forget plot] {9.7};
        \addplot[thick, dashed, domain=0.001:1000, samples=2, forget plot] {19.5};

        \node[rotate=32, above, font=\small\bfseries] at (axis cs:2.2,4.1)
        {  2039 GB/s};
        \node[above right, font=\small\bfseries] at (axis cs:1.05,9.7)
        {CUDA 9.7 TFLOP/s};
        \node[above right, font=\small\bfseries] at (axis cs:1.05,19.5)
        {Tensor 19.5 TFLOP/s};

        \addplot[black, mark=*, mark size=3pt, only marks] coordinates {(3.5,6.017)};
        \addlegendentry{Two-axis even--odd}
        \addplot[violet, mark=o, mark options={fill=none, line width=1.2pt}, mark size=3pt, only marks]
        coordinates {(5.3,9.081)};
        \addlegendentry{One-axis even--odd}
        \addplot[red, mark=square*, mark size=3pt, mark options={line width=1pt}, only marks]
        coordinates {(8.1,9.257)};
        \addlegendentry{Dense matrix--matrix}
        \addplot[orange, mark=triangle*, mark size=4pt, mark options={line width=1pt}, only marks]
        coordinates {(11.3,9.949)};
        \addlegendentry{FDM}
    \end{loglogaxis}
\end{tikzpicture}
     \caption{Roofline of the cubic three-dimensional biharmonic operator realizations.
        Arithmetic intensities and FLOP rates are NVIDIA Nsight Compute counter measurements and therefore include
        issued arithmetic that the algorithm does not use. Solid lines are the memory-bandwidth slope and the vector
        FP64 ceiling; the dashed line is the FP64 tensor-core ceiling.}
    \label{fig:biharmonic_roofline}
\end{figure}

Figure~\ref{fig:biharmonic_roofline} places the four realizations on the A100 roofline using hardware counters
collected with NVIDIA Nsight Compute. The plotted operation counts measure issued arithmetic, including basis
transformations, and differ from the dense-product counts in Table~\ref{tab:biharmonic_even_odd_variants}. The FLOP
rates therefore describe executed work rather than a common operation count for the four algorithms.

The two-axis realization has arithmetic intensity $3.5$\,FLOP/byte and reaches $6.02$\,TFLOP/s, approximately $85\%$ of
the bandwidth roof at that intensity. FDM has the highest intensity, $11.3$\,FLOP/byte, and reaches $9.95$\,TFLOP/s,
exceeding the vector FP64 ceiling because its contractions run on the tensor cores. The one-axis realization has
intensity $5.3$\,FLOP/byte and, like the two-axis one, stays below the bandwidth roof, so the additional arithmetic
does not change its run time appreciably.

\section{Smoothers}
\label{sec:biharmonic_smoothing}

Each smoother consumes a residual in the cell-wise dual representation and returns a correction that is a fixed point
of $C_\ell$, the requirement~\eqref{eq:biharmonic_preconditioner_requirement} placed on the preconditioner. Its
level-local realization is restricted to $S_\ell$ and leaves $E_\ell$ unchanged, including both value and derivative
trace layers. The diagonal, local solves, and spectral bounds must be derived for the biharmonic operator.

\subsection{Chebyshev--Jacobi Smoother}
\label{sec:biharmonic_chebyshev_jacobi}

\paragraph{Direct stiffness summation.}
For the Cartesian level grid, the map from a cell-wise dual residual to a $C^1$ primal vector is the \emph{direct
    stiffness summation} (DSS) operator~\cite{wichrowski2026DSS}
\begin{equation}
    S_\ell
    :=\mathcal{G}_\ell\mathcal{I}_{\mathrm{Riesz}}\mathcal{G}_\ell^T
    \colon V_{\mathrm{cell},\ell}^*\longrightarrow V_{C^1,\ell},
    \label{eq:biharmonic_dss}
\end{equation}
where $\mathcal{I}_{\mathrm{Riesz}}$ is the Euclidean Riesz map between assembled dual and primal spaces.
The assembly $\mathcal{G}_\ell^T$ adds the cell contributions belonging to one global DoF, and
$\mathcal{G}_\ell$ replicates the sum to every incident cell.

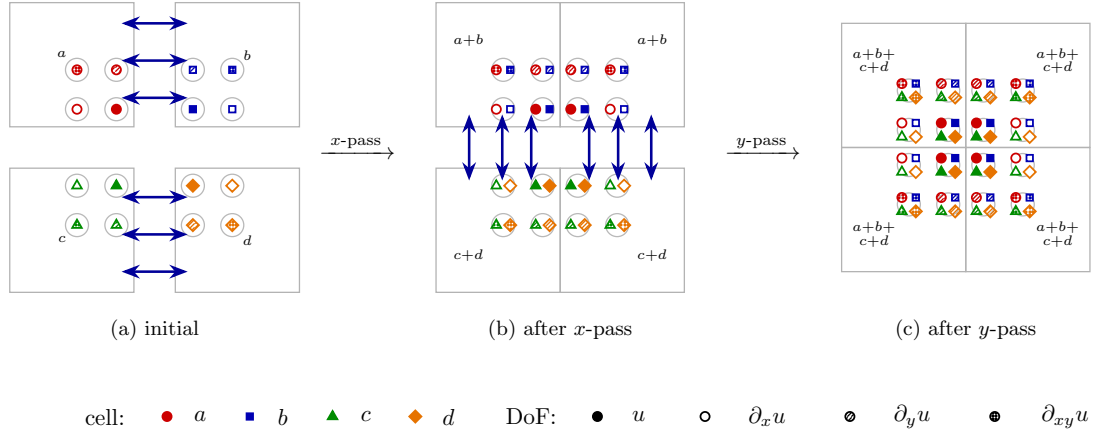
\begin{figure}[htbp]
    \centering
    \adjustbox{max width=\textwidth}{%
\usetikzlibrary{arrows.meta, positioning, calc, shapes.geometric, patterns.meta}
    \begin{tikzpicture}[
        scale=1.26,
        every node/.style={font=\small},
        gridline/.style={gray!60, line width=0.56pt},
        slot/.style={draw=gray!55, line width=0.56pt, fill=white, circle, minimum size=3.5mm, inner sep=0pt},
        exch/.style={{Stealth[length=2.8mm]}-{Stealth[length=2.8mm]}, line width=1.12pt, blue!60!black},
        ]
        \colorlet{cA}{red!80!black}
        \colorlet{cB}{blue!70!black}
        \colorlet{cC}{green!55!black}
        \colorlet{cD}{orange!90!black}
        \tikzset{
            mA/.style={circle, draw=cA, inner sep=1.4pt, line width=0.77pt},
            mB/.style={rectangle, draw=cB, inner sep=1.4pt, line width=0.77pt},
            mC/.style={regular polygon, regular polygon sides=3, draw=cC, inner sep=0.85pt, line width=0.77pt},
            mD/.style={diamond, draw=cD, inner sep=1.4pt, line width=0.77pt},
            AV/.style={fill=cA},
            AX/.style={fill=white},
            AY/.style={fill=white, pattern={Lines[angle=45,distance=1.12pt,line width=0.42pt]}, pattern color=cA},
            AXY/.style={fill=white, pattern={Dots[distance=1.26pt,radius=0.49pt]}, pattern color=cA},
            BV/.style={fill=cB},
            BX/.style={fill=white},
            BY/.style={fill=white, pattern={Lines[angle=45,distance=1.12pt,line width=0.42pt]}, pattern color=cB},
            BXY/.style={fill=white, pattern={Dots[distance=1.26pt,radius=0.49pt]}, pattern color=cB},
            CV/.style={fill=cC},
            CX/.style={fill=white},
            CY/.style={fill=white, pattern={Lines[angle=45,distance=1.12pt,line width=0.42pt]}, pattern color=cC},
            CXY/.style={fill=white, pattern={Dots[distance=1.26pt,radius=0.49pt]}, pattern color=cC},
            DV/.style={fill=cD},
            DX/.style={fill=white},
            DY/.style={fill=white, pattern={Lines[angle=45,distance=1.12pt,line width=0.42pt]}, pattern color=cD},
            DXY/.style={fill=white, pattern={Dots[distance=1.26pt,radius=0.49pt]}, pattern color=cD},
            leg/.style={circle, draw=black, inner sep=1.4pt, line width=0.77pt},
            legV/.style={fill=black},
            legX/.style={fill=white},
            legY/.style={fill=white, pattern={Lines[angle=45,distance=1.12pt,line width=0.42pt]}, pattern color=black},
            legXY/.style={fill=white, pattern={Dots[distance=1.26pt,radius=0.49pt]}, pattern color=black},
        }
        \def\bo{0.238} %
        \def\mo{0.084} %
        \newcommand{\mA}[2]{\node[mA, #1] at (#2) {};}
        \newcommand{\mB}[2]{\node[mB, #1] at (#2) {};}
        \newcommand{\mC}[2]{\node[mC, #1] at (#2) {};}
        \newcommand{\mD}[2]{\node[mD, #1] at (#2) {};}
        \newcommand{\blockslots}[1]{%
            \node[slot] at ($#1+(-\bo,\bo)$) {};
            \node[slot] at ($#1+(\bo,\bo)$) {};
            \node[slot] at ($#1+(-\bo,-\bo)$) {};
            \node[slot] at ($#1+(\bo,-\bo)$) {};}
        \newcommand{\copyAone}[1]{\blockslots{#1}
            \mA{AXY}{$#1+(-\bo,\bo)$} \mA{AY}{$#1+(\bo,\bo)$}
            \mA{AX}{$#1+(-\bo,-\bo)$} \mA{AV}{$#1+(\bo,-\bo)$};}
        \newcommand{\copyBone}[1]{\blockslots{#1}
            \mB{BY}{$#1+(-\bo,\bo)$} \mB{BXY}{$#1+(\bo,\bo)$}
            \mB{BV}{$#1+(-\bo,-\bo)$} \mB{BX}{$#1+(\bo,-\bo)$};}
        \newcommand{\copyCone}[1]{\blockslots{#1}
            \mC{CX}{$#1+(-\bo,\bo)$} \mC{CV}{$#1+(\bo,\bo)$}
            \mC{CXY}{$#1+(-\bo,-\bo)$} \mC{CY}{$#1+(\bo,-\bo)$};}
        \newcommand{\copyDone}[1]{\blockslots{#1}
            \mD{DV}{$#1+(-\bo,\bo)$} \mD{DX}{$#1+(\bo,\bo)$}
            \mD{DY}{$#1+(-\bo,-\bo)$} \mD{DXY}{$#1+(\bo,-\bo)$};}
        \newcommand{\copyAABtwo}[1]{\blockslots{#1}
            \mA{AXY}{$#1+(-\bo,\bo)+(-\mo,0)$} \mB{BXY}{$#1+(-\bo,\bo)+(\mo,0)$}
            \mA{AY}{$#1+(\bo,\bo)+(-\mo,0)$} \mB{BY}{$#1+(\bo,\bo)+(\mo,0)$}
            \mA{AX}{$#1+(-\bo,-\bo)+(-\mo,0)$} \mB{BX}{$#1+(-\bo,-\bo)+(\mo,0)$}
            \mA{AV}{$#1+(\bo,-\bo)+(-\mo,0)$} \mB{BV}{$#1+(\bo,-\bo)+(\mo,0)$};}
        \newcommand{\copyBABtwo}[1]{\blockslots{#1}
            \mA{AY}{$#1+(-\bo,\bo)+(-\mo,0)$} \mB{BY}{$#1+(-\bo,\bo)+(\mo,0)$}
            \mA{AXY}{$#1+(\bo,\bo)+(-\mo,0)$} \mB{BXY}{$#1+(\bo,\bo)+(\mo,0)$}
            \mA{AV}{$#1+(-\bo,-\bo)+(-\mo,0)$} \mB{BV}{$#1+(-\bo,-\bo)+(\mo,0)$}
            \mA{AX}{$#1+(\bo,-\bo)+(-\mo,0)$} \mB{BX}{$#1+(\bo,-\bo)+(\mo,0)$};}
        \newcommand{\copyCCDtwo}[1]{\blockslots{#1}
            \mC{CX}{$#1+(-\bo,\bo)+(-\mo,0)$} \mD{DX}{$#1+(-\bo,\bo)+(\mo,0)$}
            \mC{CV}{$#1+(\bo,\bo)+(-\mo,0)$} \mD{DV}{$#1+(\bo,\bo)+(\mo,0)$}
            \mC{CXY}{$#1+(-\bo,-\bo)+(-\mo,0)$} \mD{DXY}{$#1+(-\bo,-\bo)+(\mo,0)$}
            \mC{CY}{$#1+(\bo,-\bo)+(-\mo,0)$} \mD{DY}{$#1+(\bo,-\bo)+(\mo,0)$};}
        \newcommand{\copyDCDtwo}[1]{\blockslots{#1}
            \mC{CV}{$#1+(-\bo,\bo)+(-\mo,0)$} \mD{DV}{$#1+(-\bo,\bo)+(\mo,0)$}
            \mC{CX}{$#1+(\bo,\bo)+(-\mo,0)$} \mD{DX}{$#1+(\bo,\bo)+(\mo,0)$}
            \mC{CY}{$#1+(-\bo,-\bo)+(-\mo,0)$} \mD{DY}{$#1+(-\bo,-\bo)+(\mo,0)$}
            \mC{CXY}{$#1+(\bo,-\bo)+(-\mo,0)$} \mD{DXY}{$#1+(\bo,-\bo)+(\mo,0)$};}
        \newcommand{\copyAABCD}[1]{\blockslots{#1}
            \mA{AXY}{$#1+(-\bo,\bo)+(-\mo,\mo)$} \mB{BXY}{$#1+(-\bo,\bo)+(\mo,\mo)$}
            \mC{CXY}{$#1+(-\bo,\bo)+(-\mo,-\mo)$} \mD{DXY}{$#1+(-\bo,\bo)+(\mo,-\mo)$}
            \mA{AY}{$#1+(\bo,\bo)+(-\mo,\mo)$} \mB{BY}{$#1+(\bo,\bo)+(\mo,\mo)$}
            \mC{CY}{$#1+(\bo,\bo)+(-\mo,-\mo)$} \mD{DY}{$#1+(\bo,\bo)+(\mo,-\mo)$}
            \mA{AX}{$#1+(-\bo,-\bo)+(-\mo,\mo)$} \mB{BX}{$#1+(-\bo,-\bo)+(\mo,\mo)$}
            \mC{CX}{$#1+(-\bo,-\bo)+(-\mo,-\mo)$} \mD{DX}{$#1+(-\bo,-\bo)+(\mo,-\mo)$}
            \mA{AV}{$#1+(\bo,-\bo)+(-\mo,\mo)$} \mB{BV}{$#1+(\bo,-\bo)+(\mo,\mo)$}
            \mC{CV}{$#1+(\bo,-\bo)+(-\mo,-\mo)$} \mD{DV}{$#1+(\bo,-\bo)+(\mo,-\mo)$};}
        \newcommand{\copyBABCD}[1]{\blockslots{#1}
            \mA{AY}{$#1+(-\bo,\bo)+(-\mo,\mo)$} \mB{BY}{$#1+(-\bo,\bo)+(\mo,\mo)$}
            \mC{CY}{$#1+(-\bo,\bo)+(-\mo,-\mo)$} \mD{DY}{$#1+(-\bo,\bo)+(\mo,-\mo)$}
            \mA{AXY}{$#1+(\bo,\bo)+(-\mo,\mo)$} \mB{BXY}{$#1+(\bo,\bo)+(\mo,\mo)$}
            \mC{CXY}{$#1+(\bo,\bo)+(-\mo,-\mo)$} \mD{DXY}{$#1+(\bo,\bo)+(\mo,-\mo)$}
            \mA{AV}{$#1+(-\bo,-\bo)+(-\mo,\mo)$} \mB{BV}{$#1+(-\bo,-\bo)+(\mo,\mo)$}
            \mC{CV}{$#1+(-\bo,-\bo)+(-\mo,-\mo)$} \mD{DV}{$#1+(-\bo,-\bo)+(\mo,-\mo)$}
            \mA{AX}{$#1+(\bo,-\bo)+(-\mo,\mo)$} \mB{BX}{$#1+(\bo,-\bo)+(\mo,\mo)$}
            \mC{CX}{$#1+(\bo,-\bo)+(-\mo,-\mo)$} \mD{DX}{$#1+(\bo,-\bo)+(\mo,-\mo)$};}
        \newcommand{\copyCABCD}[1]{\blockslots{#1}
            \mA{AX}{$#1+(-\bo,\bo)+(-\mo,\mo)$} \mB{BX}{$#1+(-\bo,\bo)+(\mo,\mo)$}
            \mC{CX}{$#1+(-\bo,\bo)+(-\mo,-\mo)$} \mD{DX}{$#1+(-\bo,\bo)+(\mo,-\mo)$}
            \mA{AV}{$#1+(\bo,\bo)+(-\mo,\mo)$} \mB{BV}{$#1+(\bo,\bo)+(\mo,\mo)$}
            \mC{CV}{$#1+(\bo,\bo)+(-\mo,-\mo)$} \mD{DV}{$#1+(\bo,\bo)+(\mo,-\mo)$}
            \mA{AXY}{$#1+(-\bo,-\bo)+(-\mo,\mo)$} \mB{BXY}{$#1+(-\bo,-\bo)+(\mo,\mo)$}
            \mC{CXY}{$#1+(-\bo,-\bo)+(-\mo,-\mo)$} \mD{DXY}{$#1+(-\bo,-\bo)+(\mo,-\mo)$}
            \mA{AY}{$#1+(\bo,-\bo)+(-\mo,\mo)$} \mB{BY}{$#1+(\bo,-\bo)+(\mo,\mo)$}
            \mC{CY}{$#1+(\bo,-\bo)+(-\mo,-\mo)$} \mD{DY}{$#1+(\bo,-\bo)+(\mo,-\mo)$};}
        \newcommand{\copyDABCD}[1]{\blockslots{#1}
            \mA{AV}{$#1+(-\bo,\bo)+(-\mo,\mo)$} \mB{BV}{$#1+(-\bo,\bo)+(\mo,\mo)$}
            \mC{CV}{$#1+(-\bo,\bo)+(-\mo,-\mo)$} \mD{DV}{$#1+(-\bo,\bo)+(\mo,-\mo)$}
            \mA{AX}{$#1+(\bo,\bo)+(-\mo,\mo)$} \mB{BX}{$#1+(\bo,\bo)+(\mo,\mo)$}
            \mC{CX}{$#1+(\bo,\bo)+(-\mo,-\mo)$} \mD{DX}{$#1+(\bo,\bo)+(\mo,-\mo)$}
            \mA{AY}{$#1+(-\bo,-\bo)+(-\mo,\mo)$} \mB{BY}{$#1+(-\bo,-\bo)+(\mo,\mo)$}
            \mC{CY}{$#1+(-\bo,-\bo)+(-\mo,-\mo)$} \mD{DY}{$#1+(-\bo,-\bo)+(\mo,-\mo)$}
            \mA{AXY}{$#1+(\bo,-\bo)+(-\mo,\mo)$} \mB{BXY}{$#1+(\bo,-\bo)+(\mo,\mo)$}
            \mC{CXY}{$#1+(\bo,-\bo)+(-\mo,-\mo)$} \mD{DXY}{$#1+(\bo,-\bo)+(\mo,-\mo)$};}

        \newcommand{\lbls}[4]{%
            \node[font=\scriptsize] at ($(nw)+(-0.66,0.60)$) {$\substack{#1}$};
            \node[font=\scriptsize] at ($(ne)+(0.66,0.60)$)  {$\substack{#2}$};
            \node[font=\scriptsize] at ($(sw)+(-0.66,-0.60)$) {$\substack{#3}$};
            \node[font=\scriptsize] at ($(se)+(0.66,-0.60)$)  {$\substack{#4}$};}
        \newcommand{\lblsc}[4]{%
            \node[font=\scriptsize] at ($(nw)+(-0.42,0.42)$) {$\substack{#1}$};
            \node[font=\scriptsize] at ($(ne)+(0.42,0.42)$)  {$\substack{#2}$};
            \node[font=\scriptsize] at ($(sw)+(-0.42,-0.42)$) {$\substack{#3}$};
            \node[font=\scriptsize] at ($(se)+(0.42,-0.42)$)  {$\substack{#4}$};}

        \newcommand{\dsscells}[3]{%
            \pgfmathsetmacro\gx{#1/2}%
            \pgfmathsetmacro\gy{#2/2}%
            \begin{scope}[shift={(-\gx,\gy)}]   %
                \draw[gridline] (0,1.5) rectangle (1.5,3); \coordinate (nw) at (1.05,1.95);
            \end{scope}
            \begin{scope}[shift={(\gx,\gy)}]    %
                \draw[gridline] (1.5,1.5) rectangle (3,3); \coordinate (ne) at (1.95,1.95);
            \end{scope}
            \begin{scope}[shift={(-\gx,-\gy)}]  %
                \draw[gridline] (0,0) rectangle (1.5,1.5); \coordinate (sw) at (1.05,1.05);
            \end{scope}
            \begin{scope}[shift={(\gx,-\gy)}]   %
                \draw[gridline] (1.5,0) rectangle (3,1.5); \coordinate (se) at (1.95,1.05);
            \end{scope}
            \node[align=center, font=\footnotesize] at (1.5,-0.7) {#3};}

        \begin{scope}[shift={(0,0)}]
            \dsscells{0.5}{0.5}{(a) initial}
            \copyAone{(nw)} \copyBone{(ne)} \copyCone{(sw)} \copyDone{(se)}
            \lblsc{a}{b}{c}{d}
            \foreach \y in {2.1,2.55,3.0} \draw[exch] (1.1,\y) -- (1.9,\y);
            \foreach \y in {0.9,0.45,0.0} \draw[exch] (1.1,\y) -- (1.9,\y);
        \end{scope}

        \node at (3.95,1.5) {$\xrightarrow{\;x\text{-pass}\;}$};

        \begin{scope}[shift={(4.9,0)}]
            \dsscells{0}{0.5}{(b) after $x$-pass}
            \copyAABtwo{(nw)} \copyBABtwo{(ne)}
            \copyCCDtwo{(sw)} \copyDCDtwo{(se)}
            \lbls{a{+}b}{a{+}b}{c{+}d}{c{+}d}
            \foreach \x in {0.4,0.8,1.15} \draw[exch] (\x,1.9) -- (\x,1.1);
            \foreach \x in {1.85,2.2,2.6} \draw[exch] (\x,1.9) -- (\x,1.1);
        \end{scope}

        \node at (8.85,1.5) {$\xrightarrow{\;y\text{-pass}\;}$};

        \begin{scope}[shift={(9.8,0)}]
            \dsscells{0}{0}{(c) after $y$-pass}
            \copyAABCD{(nw)} \copyBABCD{(ne)} \copyCABCD{(sw)} \copyDABCD{(se)}
            \lbls{a{+}b{+}\\c{+}d}{a{+}b{+}\\c{+}d}{a{+}b{+}\\c{+}d}{a{+}b{+}\\c{+}d}
        \end{scope}

        \begin{scope}[shift={(6.4,-1.75)}, font=\scriptsize]
            \node at (-5.5,0) {cell:};
            \node[mA, AV] at (-4.75,0) {}; \node at (-4.35,0) {$a$};
            \node[mB, BV] at (-3.75,0) {}; \node at (-3.35,0) {$b$};
            \node[mC, CV] at (-2.75,0) {}; \node at (-2.35,0) {$c$};
            \node[mD, DV] at (-1.75,0) {}; \node at (-1.35,0) {$d$};
            \node at (-0.35,0) {DoF:};
            \node[leg, legV] at (0.45,0) {}; \node at (0.95,0) {$u$};
            \node[leg, legX] at (1.75,0) {}; \node at (2.5,0) {$\partial_x u$};
            \node[leg, legY] at (3.5,0) {}; \node at (4.25,0) {$\partial_y u$};
            \node[leg, legXY] at (5.25,0) {}; \node at (6.1,0) {$\partial_{xy} u$};
        \end{scope}

    \end{tikzpicture}%
}
    \caption{Dimensionally split DSS in two dimensions for cubic Hermite elements. At the vertex shared by four
        cells each cell carries a $2{\times}2$ block of DoFs: the value $u$, the first derivatives $\partial_x u$ and
        $\partial_y u$, and the mixed derivative $\partial_{xy} u$. Each block is drawn as seen from its own cell, so
        the four blocks are mirror images and facing slots across an interface hold the same Hermite DoF. Marker shape
        and color mark the source cell $a$--$d$; marker fill marks the Hermite DoF (solid, hollow, hatched, dotted).
        The first pass forms the horizontal pair sums and the second propagates them vertically, after which every copy
        of every DoF contains $a+b+c+d$. Each pass exchanges the two trace layers determining the pair
        $(u,\partial_n u)$, pairing value with value (layers $0$ and $p$) and derivative with derivative (layers $1$
        and $p-1$), so the vertex block is assembled without ever being addressed explicitly.}
    \label{fig:biharmonic_dss_propagation}
\end{figure}

DSS is evaluated without forming a global vector. One pass along coordinate direction $\alpha$ visits every Cartesian
interface normal to that direction, adds the two coincident copies of each interface DoF, and writes the sum to both
cells. The passes are sequential. Consequently, a later pass propagates the partial sums produced by the earlier
passes: two passes assemble a vertex shared by four cells in two dimensions
(Figure~\ref{fig:biharmonic_dss_propagation}); in three dimensions, the two directions perpendicular to an edge
assemble its four copies, and all three directions assemble the eight copies at a vertex. The exchange partners have
fixed offsets on a Cartesian grid, and each pair is owned by one thread; neither indirect addressing nor atomic
additions are required. This construction is restricted here to Cartesian level grids.

The only change from nodal $C^0$ elements is the width of the exchanged trace. A face-normal pass processes the two DoF
layers that determine the Hermite trace pair $(u,\partial_n u)$ instead of only the endpoint-value layer, pairing value
with value and derivative with derivative across the interface. Because the axis-wise cascade acts on every entry of
both layers, their edge and vertex intersections are assembled as byproducts of the same passes. Degrees of freedom on
$E_\ell$ are frozen, so DSS is needed only across matching interfaces of the level grid.

\paragraph{Jacobi stage.}
Let $A_{C^1,\ell}=\mathcal{G}_\ell^T A_\ell\mathcal{G}_\ell$, $D_{C^1,\ell}=\operatorname{diag}(A_{C^1,\ell})$, and let
$D_{\mathrm{cell},\ell}$ contain a copy of this diagonal on every incident cell. It is constructed once per level by
placing the scaled cell diagonals in cell-wise storage and summing their coincident entries with DSS. Let $\kappa$ be
the number of cell copies of each free level DoF, $K_{\mathrm{cell},\ell}$ its replicated diagonal, and $M_{S_\ell}$
the mask that vanishes on $E_\ell$. The stored inverse diagonal is $\widehat
    D_\ell^{-1}=M_{S_\ell}K_{\mathrm{cell},\ell}D_{\mathrm{cell},\ell}^{-1}M_{S_\ell}$; its zero entries freeze both
constrained trace layers.

Define averaging DSS on the matching level interfaces by
\begin{equation}
    \bar S_\ell
    :=\mathcal{G}_\ell(\mathcal{G}_\ell^T\mathcal{G}_\ell)^{-1}\mathcal{G}_\ell^T,
    \qquad
    \bar S_\ell K_{\mathrm{cell},\ell}D_{\mathrm{cell},\ell}^{-1}
    =D_{\mathrm{cell},\ell}^{-1}S_\ell.
    \label{eq:biharmonic_averaging_dss}
\end{equation}
Starting from a continuous iterate, one damped Jacobi stage is
\begin{equation}
    \tilde r_\ell^{(j-1)}=\tilde b_\ell-A_\ell u_\ell^{(j-1)},
    \qquad
    u_\ell^{(j)}
    =\bar S_\ell\!\left(u_\ell^{(j-1)}
    +\omega_j\widehat D_\ell^{-1}\tilde r_\ell^{(j-1)}\right).
    \label{eq:biharmonic_jacobi_stage}
\end{equation}
where the update inside parentheses uses the unassembled cell residual. By~\eqref{eq:biharmonic_averaging_dss}, the
final averaging pass both assembles this update and restores $C^1$ continuity. A stage therefore comprises the
cell-wise operator application, a pointwise Jacobi update, and one two-layer DSS operation.

\paragraph{Chebyshev acceleration.}
Chebyshev acceleration replaces a fixed damping parameter by $n$ Richardson weights for a target interval
$[a,b]\subset(0,\infty)$ of the spectrum of the diagonally preconditioned operator on $S_\ell$. With
\begin{equation}
    c=\frac{a+b}{2},\qquad \delta=\frac{b-a}{2},\qquad
    \omega_j=\frac{1}{c-\delta\cos\!\left(\frac{(2j-1)\pi}{2n}\right)},
    \quad j=1,\ldots,n,
    \label{eq:biharmonic_chebyshev_weights}
\end{equation}
the stages~\eqref{eq:biharmonic_jacobi_stage} give the error polynomial
\begin{equation}
    p_n(\lambda)=\prod_{j=1}^n(1-\omega_j\lambda)
    =\frac{T_n((c-\lambda)/\delta)}{T_n(c/\delta)},
    \qquad
    \max_{\lambda\in[a,b]}|p_n(\lambda)|=|T_n(c/\delta)|^{-1},
    \label{eq:biharmonic_chebyshev_polynomial}
\end{equation}
where $T_n$ is the Chebyshev polynomial of the first kind. The factorized form changes only the scalar weight between
Jacobi stages and needs no auxiliary recurrence vector.

The spectral interval is fixed rather than estimated during setup. We use $[a,b]=[0.2,2]$: the upper endpoint safely
covers the observed value $\lambda_{\max}(D_{C^1,\ell}^{-1}A_{C^1,\ell})\approx1.96$, while the lower endpoint selects
the upper decade of the spectrum for smoothing and leaves lower-frequency components to the coarse-grid correction. The
reported Chebyshev--Jacobi experiments use $n=3$ stages in each pre- and post-smoothing application; the coarsest level
uses the same construction with $20$ stages in place of an exact solve.

\subsection{Vertex-Patch Smoother}
\label{sec:biharmonic_patch_smoother}

For an interior level vertex $v$, let
\begin{equation}
    \omega_v:=\{K\in\mathcal{T}_\ell:v\in\overline K\}.
    \label{eq:biharmonic_vertex_patch}
\end{equation}
On a Cartesian mesh, $\omega_v$ contains the eight cells meeting at $v$. A cubic tensor-product Hermite element
has two DoFs at each endpoint in each coordinate direction: a value and a derivative. Their tensor
product gives $2^3=8$ DoFs at the shared vertex. These eight DoFs form the vertex block;
all other DoFs in the patch are held fixed during its update.

The patch matrix $A_\ell^{\mathrm{patch}}\in\mathbb{R}^{8\times8}$ is obtained by restricting each of the eight cell
matrices to the DoFs at the shared vertex and summing their contributions. This is the assembled biharmonic operator on
the vertex block with a zero correction prescribed on all remaining patch DoFs. On a uniform level, every interior
vertex has the same cell configuration and cell size. The same matrix and its dense $8\times8$ inverse can therefore be
used for every patch on that level.

Patch setup is performed once per level. The cell--vertex incidences are sorted by vertex, vertices incident to eight
owned level cells are retained, and the eight incident cells are stored in a fixed corner order. Patches are then
grouped into eight colors, and $(A_\ell^{\mathrm{patch}})^{-1}$ is computed. At run time, a patch is represented only
by its eight cell indices; the corner position identifies the value and derivative entries belonging to the shared
vertex.

For every cell $K\in\omega_v$, the smoother evaluates $A_Ku_K$ and subtracts it from the cell-local right-hand side. It
then extracts the eight residual entries belonging to $v$ and sums their copies over the eight cells. Denoting this
restriction by $|_v$, the condensed patch residual and its correction are
\begin{equation}
    r_v:=\sum_{K\in\omega_v}
    \left.\bigl(\widetilde f_K-A_Ku_K\bigr)\right|_v,
    \qquad
    \delta_v:=\omega (A_\ell^{\mathrm{patch}})^{-1}r_v.
    \label{eq:biharmonic_patch_correction}
\end{equation}
Here $r_v,\delta_v\in\mathbb{R}^8$, $\omega$ is a damping parameter, and $\omega=1$ gives the undamped patch solve.
The correction is scattered to the eight cells by adding it to every cell-local copy of the vertex block,
\begin{equation}
    \left.u_K\right|_v\leftarrow\left.u_K\right|_v+\delta_v,
    \qquad K\in\omega_v,
    \label{eq:biharmonic_patch_scatter}
\end{equation}
so all copies receive the same update. Consequently, an initially $C^1$ iterate remains $C^1$ without a separate
continuity-restoration operation.

The patches are partitioned into eight colors according to the parity of their vertex coordinates on the dyadic grid.
Patches of one color share no cell and can therefore be processed concurrently. Colors are visited successively, so the
residual for one color includes all corrections made by earlier colors; a full pass is thus a multiplicative
vertex-block sweep. Reversing the color order gives its adjoint. Forward pre-smoothing followed by reverse
post-smoothing yields the symmetric smoothing composition required when the V-cycle is used as a conjugate-gradient
preconditioner.

Because no cell outside $\omega_v$ contains a copy of a DoF owned by $v$, the residual of that block immediately after
its update is
\begin{equation}
    r_v^{\mathrm{new}}
    =r_v-A_\ell^{\mathrm{patch}}\delta_v
    =(1-\omega)r_v.
    \label{eq:biharmonic_patch_residual_after_update}
\end{equation}
Thus the undamped update annihilates the assembled residual of the current vertex block. Corrections made by subsequent
colors can change it again, giving the multiplicative coupling between the vertex blocks.

Only vertices whose complete block belongs to $S_\ell$ are admissible patch centers. In particular, a boundary vertex
is omitted because the clamped conditions fix its value and normal-derivative layers. At a refinement edge, omitting
every patch whose vertex block intersects $E_\ell$ prevents the smoother from changing dependent value or derivative
DoFs. In both cases the rule reduces to the same test, selecting the vertices with all eight incident cells owned by
the level, the test the setup above performs. A level of a locally refined hierarchy is still a uniform cut of the
domain, so it needs no patch matrices beyond the one it already builds; the rule simply retains fewer patches on the
levels that carry a refinement edge.

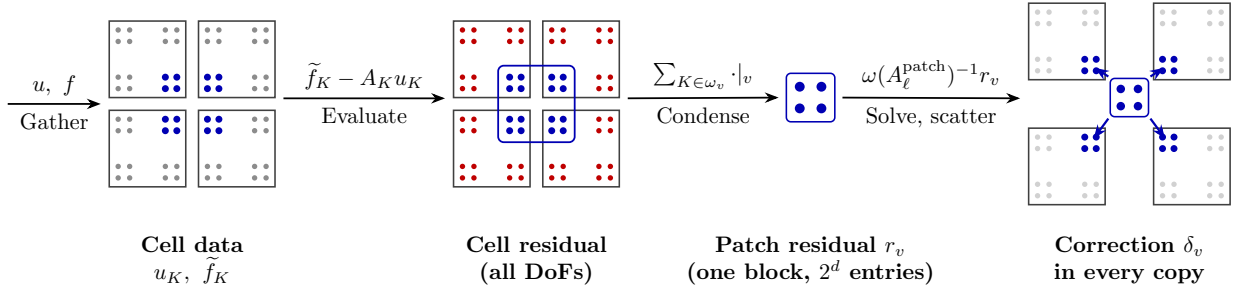
\begin{figure}[htbp]
    \centering
    \resizebox{\textwidth}{!}{%
\usetikzlibrary{arrows.meta, positioning, calc}
\begin{tikzpicture}[
        scale=1.0,
        every node/.style={font=\small},
        cell/.style={draw=black!75, line width=0.7pt, minimum size=1.2cm, inner sep=0pt},
        dof/.style={circle, inner sep=0pt, minimum size=2.4pt},
        frozen/.style={dof, fill=black!45},
        resid/.style={dof, fill=red!75!black},
        block/.style={dof, minimum size=3.1pt, fill=blue!70!black},
        faint/.style={dof, fill=black!18},
        arrow/.style={-{Stealth[length=2.2mm]}, line width=0.8pt, shorten >=1.2mm, shorten <=1.2mm},
        lab/.style={font=\small\bfseries, align=center},
    ]
    \newcommand{\cellDofs}[2]{%
        \foreach \sx in {-1,1}{\foreach \sy in {-1,1}{%
                        \foreach \ux in {0.46,0.28}{\foreach \uy in {0.46,0.28}{%
                                        \node[#2] at ($(#1.center)+(\sx*\ux,\sy*\uy)$) {};}}}}}
    \newcommand{\cornerBlock}[4]{%
        \foreach \ux in {0.46,0.28}{\foreach \uy in {0.46,0.28}{%
                        \node[#4] at ($(#1.center)+(#2*\ux,#3*\uy)$) {};}}}
    \newcommand{\panel}[4]{%
        \node[cell] (#1a) at ($(#2)+(-#4,#4)$) {};
        \node[cell] (#1b) at ($(#2)+(#4,#4)$) {};
        \node[cell] (#1c) at ($(#2)+(-#4,-#4)$) {};
        \node[cell] (#1d) at ($(#2)+(#4,-#4)$) {};
        \foreach \c in {a,b,c,d}{\cellDofs{#1\c}{#3}}}

    \panel{pA}{0,0}{frozen}{0.70}
    \cornerBlock{pAa}{1}{-1}{block}\cornerBlock{pAb}{-1}{-1}{block}
    \cornerBlock{pAc}{1}{1}{block}\cornerBlock{pAd}{-1}{1}{block}
    \node[lab, anchor=north] at (0,-1.95) {Cell data\\$u_K,\ \widetilde f_K$};

    \panel{pB}{5.4,0}{resid}{0.70}
    \cornerBlock{pBa}{1}{-1}{block}\cornerBlock{pBb}{-1}{-1}{block}
    \cornerBlock{pBc}{1}{1}{block}\cornerBlock{pBd}{-1}{1}{block}
    \draw[draw=blue!70!black, line width=0.8pt, rounded corners=2.5pt]
    ($(5.4,0)+(-0.60,-0.60)$) rectangle ($(5.4,0)+(0.60,0.60)$);
    \node[lab, anchor=north] at (5.4,-1.95) {Cell residual\\(all DoFs)};

    \coordinate (pC) at (9.7,0.10);
    \draw[draw=blue!70!black, line width=0.7pt, rounded corners=2pt]
    ($(pC)+(-0.38,-0.38)$) rectangle ($(pC)+(0.38,0.38)$);
    \foreach \sx in {-1,1}{\foreach \sy in {-1,1}{%
                    \node[block, minimum size=4.4pt] at ($(pC)+(\sx*0.17,\sy*0.17)$) {};}}
    \node[lab, anchor=north] at (9.7,-1.95) {Patch residual $r_v$\\(one block, $2^d$ entries)};

    \panel{pD}{14.7,0}{faint}{0.98}
    \cornerBlock{pDa}{1}{-1}{block}\cornerBlock{pDb}{-1}{-1}{block}
    \cornerBlock{pDc}{1}{1}{block}\cornerBlock{pDd}{-1}{1}{block}
    \coordinate (pD) at (14.7,0.10);
    \fill[white] ($(pD)+(-0.30,-0.30)$) rectangle ($(pD)+(0.30,0.30)$);
    \draw[draw=blue!70!black, line width=0.7pt, rounded corners=2pt]
    ($(pD)+(-0.30,-0.30)$) rectangle ($(pD)+(0.30,0.30)$);
    \foreach \sx in {-1,1}{\foreach \sy in {-1,1}{%
                    \node[block, minimum size=3.8pt] at ($(pD)+(\sx*0.135,\sy*0.135)$) {};}}
    \foreach \c/\sx/\sy in {a/1/-1,b/-1/-1,c/1/1,d/-1/1}{%
            \draw[arrow, blue!70!black, line width=1.0pt, shorten <=0pt, shorten >=0.8mm]
            ($(pD)+(-\sx*0.33,-\sy*0.33)$) -- ($(pD\c.center)+(\sx*0.37,\sy*0.37)$);}
    \node[lab, anchor=north] at (14.7,-1.95) {Correction $\delta_v$\\in every copy};

    \draw[arrow] ($(pAa.west)+(-1.7,-0.70)$) -- node[above, midway] {$u,\ f$}
    node[below, midway] {Gather} ($(pAa.west)+(0,-0.70)$);
    \draw[arrow] (pAb.east|-pAa.south) -- node[above, midway] {$\widetilde f_K-A_Ku_K$}
    node[below, midway] {Evaluate} (pBa.west|-pAa.south);
    \draw[arrow] (pBb.east|-pAa.south) -- node[above, midway] {$\sum_{K\in\omega_v}\cdot|_v$}
    node[below, midway] {Condense} ($(pC)+(-0.38,0)$);
    \draw[arrow] ($(pC)+(0.38,0)$) -- node[above, midway] {$\omega(A_\ell^{\mathrm{patch}})^{-1}r_v$}
    node[below, midway] {Solve, scatter} (pDa.west|-pAa.south);
\end{tikzpicture}
}
    \caption{Stages of one vertex-patch update, drawn in two dimensions. Each cell carries its own copy of the DoFs at the shared vertex, and
        every cell corner carries a $2^d$ block of Hermite DoFs. The
        cell residual is evaluated on all cell DoFs, the copies belonging to the patch vertex are summed into the single
        block $r_v$, which is not in the cell-wise format, the dense inverse of $A_\ell^{\mathrm{patch}}$ is
        applied, and the correction is added to every cell-local copy of the vertex block, which keeps the iterate
        $C^1$.}
    \label{fig:biharmonic_vertex_patch_stages}
\end{figure}

Algorithm~\ref{alg:biharmonic_vertex_patch_sweep} gives one sweep. All patches of the current color are independent and
are processed concurrently, whereas the colors remain sequential because each color reads the corrections written by
its predecessors.
\begin{algorithm}[H]
    \caption{Fused multiplicative vertex-patch sweep}
    \label{alg:biharmonic_vertex_patch_sweep}
    \begin{algorithmic}[1]
        \Require Primal cell-wise iterate $u_\ell$; unassembled dual right-hand side $\widetilde f_\ell$;
        color sets $\{\mathcal{V}_{\ell,c}\}_{c=0}^{7}$; patch inverse $(A_\ell^{\mathrm{patch}})^{-1}$; damping $\omega$;
        color order $\pi$.
        \Ensure Updated $C^1$ iterate $u_\ell$.
        \For{$q=0,\ldots,7$}
        \State $c\gets\pi(q)$
        \State Process all $v\in\mathcal{V}_{\ell,c}$ concurrently \Comment{One device launch}
        \ForAll{$v\in\mathcal{V}_{\ell,c}$ \textbf{ in parallel}}
        \State $r_v\gets0$
        \ForAll{$K\in\omega_v$}
        \State Load $u_K$ and $\widetilde f_K$
        \State $y_K\gets A_Ku_K$ \Comment{Cell-local operator evaluation}
        \State $r_v\gets r_v+\left.(\widetilde f_K-y_K)\right|_v$ \Comment{Condense to eight entries}
        \EndFor
        \State $\delta_v\gets\omega(A_\ell^{\mathrm{patch}})^{-1}r_v$ \Comment{Dense $8\times8$ product}
        \ForAll{$K\in\omega_v$}
        \State $\left.u_K\right|_v\gets\left.u_K\right|_v+\delta_v$ \Comment{Scatter to every copy}
        \EndFor
        \EndFor
        \EndFor
        \State \Return $u_\ell$
    \end{algorithmic}
\end{algorithm}

Within a color launch, the cell results $y_K$ and the patch residual $r_v$ in
Algorithm~\ref{alg:biharmonic_vertex_patch_sweep} remain local to the executing thread block. The eight cell
evaluations, restriction and reduction, $8\times8$ solve, and replicated scatter are fused; neither $y_K$ nor $r_v$ is
written to global memory. The cell evaluation uses the two-axis even--odd realization. This kernel is the one place in this work that requires each cell's $4\times4\times4$ DoF
cube to be contiguous in memory, since a patch gathers eight complete cubes at once; everything else, the level
operators and the assembled diagonal included, is stored in the blocked layout of~\cite{wichrowski2026DSS} with block
size $B=32$, required by the DSS passes and flat in the volumetric rate.

The fused algorithm still evaluates eight complete cell operators per patch. For each cell visit in FP64 it reads the
$64$ entries of $u_K$ and the $64$ entries of $\widetilde f_K$, and it writes only the eight entries owned by the patch
vertex. Excluding the small read-only operator and patch-inverse tables, the minimum field-data volume is therefore
\begin{equation}
    64\cdot8+64\cdot8+8\cdot8=1088\ \text{bytes per cell visit},
    \qquad
    8\cdot1088=8704\ \text{bytes per patch}.
    \label{eq:biharmonic_patch_bytes}
\end{equation}
For a uniform $n\times n\times n$ mesh, there are $(n-1)^3$ interior vertex patches and
$8(n-1)^3$ cell visits per sweep. The minimum traffic normalized by the $n^3$ mesh cells is consequently
\begin{equation}
    b_{\mathrm{sweep}}(n)
    =8704\left(1-\frac{1}{n}\right)^3
    \quad\text{bytes per mesh cell}.
    \label{eq:biharmonic_patch_bytes_per_cell}
\end{equation}
Equivalently, a mesh cell is revisited $8(1-1/n)^3$ times on average. For $n=128$, this gives $8501$ bytes per mesh
cell, compared with $1024$ bytes for one volumetric operator evaluation. The factor of $8.30$ is intrinsic to the
colored multiplicative algorithm and is separated from the measured hardware efficiency of its realization.

\paragraph{Possible extension to unstructured meshes.}
The continuity-preserving patch update suggests a generalization to unstructured meshes. On a vertex patch, one could
construct a local correction space whose functions have matching value and physical normal-derivative traces across
internal facets and vanishing traces of both quantities on the patch boundary. Extending such corrections by zero would
yield globally $C^1$ functions, so local patch updates would preserve an initially $H^2$-conforming field independently
of the mesh connectivity. This requires local finite element spaces rich enough to support nontrivial conforming patch
corrections on the chosen element geometries. Their construction, and whether the resulting collection of patch spaces
provides effective smoothing, would need to be established for the unstructured discretization.

\subsection{Performance of the Smoothers}
\label{sec:biharmonic_smoother_performance}

Table~\ref{tab:biharmonic_smoother_rates} compares throughput per smoothing application. A Chebyshev--Jacobi stage
consists of one volumetric evaluation, a pointwise update, and $d$ DSS passes. Each pass exchanges two trace layers per
participating face. An interior cell participates at both faces normal to the pass direction, covering $4/(p+1)$ of its
entries for $p\geq3$, including all entries for cubic elements. A colored patch sweep makes $8(1-1/n)^3$ cell visits
per mesh cell on average, or $7.81$ for $n=128$. Its field-data volume is $8.30$ times that of one volumetric
evaluation. \secref{sec:biharmonic_results_solve} compares total solve times.

A stage costs $6.7$ operator evaluations, so the operator is the smaller part of it. The reason is the width of the DSS
cascade at $p=3$: the two trace layers exchanged at the two faces normal to a pass direction are all four layers the
cell has in that direction, so each of the $d$ passes reads and writes the entire level vector, moving the $1024$ bytes
per cell that one operator evaluation moves. Three passes, the operator, and a pointwise update over two input vectors
account for the measured factor. This cost is fixed by the element and the mesh dimension. It is independent of how the
operator is evaluated, and it does not fall when the operator is made faster: the $1.96\times$ separating the two-axis
realization from fast diagonalization in Table~\ref{tab:biharmonic_operator_performance} acts on roughly one seventh of
a stage. The corresponding statement for the patch smoother is~\eqref{eq:biharmonic_patch_bytes}, where the revisits
rather than the synchronization set the cost.

\begin{table}[htbp]
    \centering
\begin{tabular}{@{}l r r r@{}}
    \toprule
    Smoothing application & time [ms] & [Mcells/s] & cost in $Au$ \\
    \midrule
    Chebyshev--Jacobi stage & 8.384  & 250.1 & 6.7 \\
    Vertex-patch sweep      & 14.741 & 142.3 & 11.9 \\
    \midrule
    Two-axis operator $Au$  & 1.243  & 1687.2 & 1.0 \\
    \bottomrule
\end{tabular}
     \caption{Cost of one smoothing application on a mesh of $128^3$ cells. The rate counts
        cells of the mesh, so it is the speed at which one application sweeps the mesh whatever it does inside, and the
        last column expresses the same thing as the number of volumetric operator evaluations $Au$ that fit in the time
        one application takes. The two smoothers run in different memory layouts, the blocked one for the Chebyshev
        stage and cell-contiguous storage for the patch sweep, so they are comparable per application of a smoother and
        not per byte.}
    \label{tab:biharmonic_smoother_rates}
\end{table}

The cost of the sweep is dominated by the revisits, each moving an intrinsic $1088$
bytes~\eqref{eq:biharmonic_patch_bytes}. Table~\ref{tab:biharmonic_patch_bandwidth} reports its rate in both units,
cells of the mesh per second and bytes per second, with the volumetric operator as reference.

\begin{table}[htbp]
    \centering
\begin{tabular}{@{}l r r r r@{}}
    \toprule
    Kernel & visits/cell & \multicolumn{2}{c}{rate} & of peak \\
    \cmidrule(lr){3-4}
           &             & [Mcells/s] & [GB/s]      &         \\
    \midrule
    Vertex-patch sweep           & 7.81 & 142.3  & 1209 & 59.3\% \\
    Addressing control           & 7.81 & 146.8  & 1763 & 86.4\% \\
    \midrule
    Two-axis volumetric operator & 1.00 & 1687.2 & 1728 & 84.7\% \\
    \bottomrule
\end{tabular}
     \caption{Rates of the vertex-patch sweep on a mesh of $128^3$ cells. Visits per cell is how
        often one application touches each cell of the mesh, and the cell rate counts cells of the mesh rather than
        visits. The byte rate is the intrinsic field traffic over the measured kernel time, $1088$ bytes per cell visit
        for the sweep, $1536$ for the control, and $1024$ per cell for the operator; the last column gives it as a
        fraction of peak bandwidth. The addressing control replaces the arithmetic of the sweep by $u+\omega f$
        on the same patches and tile, and stores the whole cube rather than the eight vertex entries.}
    \label{tab:biharmonic_patch_bandwidth}
\end{table}

The sweep reaches an effective bandwidth of $1209$\,GB/s, $70\%$ of the $1728$\,GB/s attained by the volumetric
operator and $59.3\%$ of peak. It requests $8.30$ times the operator's field-data volume and takes $11.9$ times as
long. These bandwidths use the minimum field-data volume; they do not measure total DRAM traffic, which may include
additional transactions for partial stores.

The sweep stores eight vertex DoFs in four $32$-byte sectors, using $64$ of the $128$ bytes covered by those sectors.
The addressing control visits the same patches but replaces the arithmetic by $u+\omega f$ and stores the full $64$-DoF
cube. It reaches an effective bandwidth of $1763$\,GB/s, showing that patch addressing is compatible with high
throughput when combined with full-cube stores. Its runtime is $14.281$\,ms, compared with $14.741$\,ms for the patch
sweep. Because both arithmetic and store volume change, this comparison does not isolate the cost of the arithmetic or
establish that partial stores account for the entire bandwidth difference.

\section{Numerical Results}
\label{sec:biharmonic_results}

We evaluate discretization error, continuity of the value and normal-derivative traces, iteration counts, and solve
times on uniform and locally refined meshes. All experiments use cubic Hermite elements, $p=3$, in three dimensions and
in FP64 arithmetic, on the device of \secref{sec:biharmonic_operator_comparison}.

\subsection{Experimental Settings}
\label{sec:biharmonic_results_settings}

The Chebyshev--Jacobi V-cycle applies one pre- and one post-smoothing application of the three-stage smoother of
\secref{sec:biharmonic_chebyshev_jacobi} on the spectral interval $[a,b]=[0.2,2]$, and replaces the coarse-level solve
by the same construction with $20$ stages. It preconditions a conjugate gradient iteration. The vertex-patch V-cycle
applies two pre-smoothing and two post-smoothing multiplicative sweeps of \secref{sec:biharmonic_patch_smoother} with
$\omega=1$, both in forward color order, and is iterated on its own. Reversing the color order of the post-smoother
gives a symmetric cycle suitable for conjugate gradient preconditioning. We report this symmetric variant both as a
stationary iteration and as a conjugate gradient preconditioner. All solver configurations stop when the pairing
$\langle\tilde r,z\rangle^{1/2}$, the measure required by the cell-wise dual representation, has decreased by a factor
$10^{-8}$. Both cycles evaluate the operator with the two-axis even--odd realization of
\secref{sec:biharmonic_even_odd}, on the blocked layout of~\cite{wichrowski2026DSS} with $B=32$; the patch kernel alone
runs on cell-contiguous storage, for the reason given in \secref{sec:biharmonic_patch_smoother}.

Throughput is $N_{\mathrm{cells}}/t$, where $t$ is the wall time of the measured operation. This normalization reflects
the cell-wise storage of $(p+1)^d=64$ entries per cell. It is also the unit in which the two discretizations compared
below are commensurable: both approximate at order $p+1$, so their interpolation error is governed by the mesh, and a
cell is the unit of accuracy delivered. The number of independent DoFs per cell is not: cubic Hermite carries $2^d=8$
of them, a $C^0$ Lagrange space of the same degree carries $p^d=27$. We therefore report both normalizations wherever
the comparison is made. Solve throughput is $N_{\mathrm{cells}}/t_{\mathrm{solve}}$ and is reported together with the
iteration count.

\subsection{Discretization and Constraint Verification}
\label{sec:biharmonic_verification}

We solve the clamped biharmonic problem on the unit cube with the manufactured solution
$u=\prod_{i=1}^{d}\sin^2(\pi x_i)$, which satisfies the clamped boundary conditions, and uniform refinement.
Table~\ref{tab:biharmonic_convergence} reports the $L^2$ error. Its measured rate is $3.97$ to $4.00$ over the
sequence, consistent with the optimal order $p+1$ for cubic elements. The maximum jumps of the trace and of the
physical normal derivative across matching interfaces are zero on every mesh in the sequence. These measurements agree
with the continuity conditions imposed by the constraint operator of \secref{sec:biharmonic_cellwise}.

Iteration counts do not grow over the tested uniform meshes. Conjugate gradients preconditioned by the
Chebyshev--Jacobi V-cycle require five or six iterations as the mesh grows from $512$ to $2\,097\,152$ cells. The
stationary patch cycle requires seven or eight iterations in its nonsymmetric form, with a contraction factor $\rho$
between $0.065$ and $0.100$, and ten in its symmetric form, at $\rho$ between $0.137$ and $0.153$. On the finest mesh,
reversing the post-smoothing color order increases $\rho$ from $0.066$ to $0.142$. Conjugate gradient acceleration
reduces the symmetric variant to six or seven iterations, with an effective contraction factor between $0.042$ and
$0.070$.

\begin{table}[htbp]
    \centering
    \adjustbox{max width=\textwidth}{%
\small
\begin{tabular}{@{}c|rr|rr|c|cc|cc|cc@{}}
    \toprule
             &         &       &                     &        & {Cheb.}  & \multicolumn{2}{c|}{Patch, nonsym.} & \multicolumn{2}{c|}{Patch, sym.} & \multicolumn{2}{c}{Patch, sym.\ + CG}                                             \\
    {$\ell$} & {Cells} & {$h$} & {$\|u-u_h\|_{L^2}$} & {rate} & {CG it.} & {V-cyc.}                            & {$\rho$}                         & {V-cyc.}                              & {$\rho$} & {CG it.} & {$\rho$} \\
    \midrule
    3 & 512       & $1.25\cdot10^{-1}$ & $1.245\cdot10^{-4}$ & ---  & 5 & 7 & 0.068 & 10 & 0.137 & 6 & 0.050 \\
    4 & 4,096     & $6.25\cdot10^{-2}$ & $7.862\cdot10^{-6}$ & 3.98 & 6 & 8 & 0.100 & 10 & 0.153 & 7 & 0.070 \\
    5 & 32,768    & $3.12\cdot10^{-2}$ & $4.924\cdot10^{-7}$ & 4.00 & 6 & 7 & 0.068 & 10 & 0.141 & 6 & 0.048 \\
    6 & 262,144   & $1.56\cdot10^{-2}$ & $3.078\cdot10^{-8}$ & 4.00 & 5 & 7 & 0.065 & 10 & 0.142 & 6 & 0.044 \\
    7 & 2,097,152 & $7.81\cdot10^{-3}$ & $1.969\cdot10^{-9}$ & 3.97 & 5 & 7 & 0.066 & 10 & 0.142 & 6 & 0.042 \\
    \bottomrule
\end{tabular}
}
    \caption{Uniform refinement of the unit cube with cubic Hermite elements in three dimensions. The error and its rate
        are those of the Chebyshev--Jacobi solve. \emph{CG it.} counts conjugate gradient iterations preconditioned by
        the Chebyshev--Jacobi V-cycle, \emph{V-cyc.} counts stationary vertex-patch cycles in the nonsymmetric and the
        symmetric color order, and the last pair counts conjugate gradient iterations preconditioned by the symmetric
        cycle. In each case $\rho$ is the measured contraction factor of the iteration reported beside it. All solvers reduce
        $\langle\tilde r,z\rangle^{1/2}$ by $10^{-8}$.}
    \label{tab:biharmonic_convergence}
\end{table}

\subsection{Solve Throughput}
\label{sec:biharmonic_results_solve}

Table~\ref{tab:biharmonic_solve} reports solve throughput. Up to $\ell=5$, throughput grows by a factor of roughly six
per refinement while the cell count grows by eight, consistent with a substantial latency contribution on small meshes.
The Chebyshev--Jacobi solver reaches $4.23$\,Mcells/s on the finest mesh, the nonsymmetric patch cycle
$3.29$\,Mcells/s, its Krylov-accelerated symmetric counterpart $3.16$\,Mcells/s, and the symmetric cycle iterated on
its own $2.30$\,Mcells/s.

On the finest mesh, reversing the color order leaves the measured cycle time nearly unchanged. Conjugate gradient
acceleration raises the cost of a symmetric patch iteration by about $20\%$ and reduces the iteration count from ten to
six. The Chebyshev--Jacobi preconditioned iteration is about $10\%$ cheaper than that accelerated patch iteration and
requires five iterations. Both iteration cost and iteration count therefore contribute to its shorter solve time.

\begin{table}[htbp]
    \centering
    \adjustbox{max width=\textwidth}{%
\begin{tabular}{@{}c|r|rr|rr|rr|rr@{}}
    \toprule
             &         & \multicolumn{2}{c|}{Chebyshev--Jacobi} & \multicolumn{2}{c|}{Patch, nonsym.} & \multicolumn{2}{c|}{Patch, sym.} & \multicolumn{2}{c}{Patch, sym.\ + CG}                                                                   \\
    {$\ell$} & {Cells} & {CG it.}                               & {[Mcells/s]}                        & {V-cyc.}                         & {[Mcells/s]}                         & {V-cyc.} & {[Mcells/s]} & {CG it.} & {[Mcells/s]} \\
    \midrule
    3 & 512       & 5 & 0.016 & 7 & 0.019 & 10 & 0.014 & 6 & 0.019 \\
    4 & 4,096     & 6 & 0.091 & 8 & 0.103 & 10 & 0.083 & 7 & 0.104 \\
    5 & 32,768    & 6 & 0.598 & 7 & 0.699 & 10 & 0.486 & 6 & 0.699 \\
    6 & 262,144   & 5 & 3.25  & 7 & 2.54  & 10 & 1.79  & 6 & 2.50  \\
    7 & 2,097,152 & 5 & 4.23  & 7 & 3.29  & 10 & 2.30  & 6 & 3.16  \\
    \bottomrule
\end{tabular}
}
    \caption{Time to solution under uniform refinement. Rates are
        $N_{\mathrm{cells}}/t_{\mathrm{solve}}$. The nonsymmetric patch cycle post-smooths in the forward color order
        and is iterated on its own; the symmetric cycle reverses that order, making it its own adjoint, and is reported
        both as a stationary iteration and as a conjugate gradient preconditioner. All solvers reduce
        $\langle\tilde r,z\rangle^{1/2}$ by $10^{-8}$.}
    \label{tab:biharmonic_solve}
\end{table}
\FloatBarrier

The patch smoother requires more data movement per application: a sweep moves $8.3$ times the data of a volumetric
operator evaluation, as derived in \secref{sec:biharmonic_smoother_performance}. A Chebyshev stage consists of an
operator evaluation, a pointwise update, and a two-layer DSS operation. Even with conjugate gradient acceleration, the
patch solver has approximately $25\%$ lower throughput than the Chebyshev--Jacobi solver on the finest mesh. This
comparison applies to the measured FP64 uniform-mesh configurations, with eight sequential color launches per patch
sweep. The patch smoother does not require a prescribed spectral interval.

For reference, the tensor-product patch multigrid solver for biharmonic problems of Witte et al.~\cite{witte2025tensor}
reports $6.4$ and $10.6$\,MDoF/s in double and in mixed precision for a $C^0$ interior penalty discretization of degree
three. At $27$ independent DoFs per cell, these rates are approximately $0.23$\,Mcells/s on a $128^3$ mesh in double
precision and $0.39$\,Mcells/s on a $256^3$ mesh in mixed precision. The $4.23$\,Mcells/s measured here is
approximately $18$ times the double-precision reference rate for the same $128^3$ mesh and polynomial degree. Per
independent DoF the same measurement is $33.8$\,MDoF/s, approximately $5.3$ times the reference. The larger ratio
follows from the element rather than from the implementation: cubic Hermite reaches the same order of approximation
with $8$ unknowns per cell instead of $27$. Both studies specify a relative stopping threshold of $10^{-8}$, but
equivalence of their stopping norms has not been established. Both use A100 GPUs with $80$\,GB of memory; the reference
reports different peak compute rates, and the reference rates are estimated from a published figure. These differences
limit the interpretation of either ratio.

With two pre-smoothing and two post-smoothing sweeps, the multiplicative vertex-patch method of Witte et
al.~\cite{witte2025tensor} requires $3.8$--$4.0$ iterations for cubic elements in three dimensions. Our symmetric
vertex-patch V-cycle uses the same number of sweeps and requires six or seven CG iterations on the tested uniform
meshes. The reference reports fractional counts derived from residual reduction, whereas we report completed
iterations.

The smaller correction spaces in our smoother provide a possible explanation for its higher iteration count. Each local
solve updates only the eight DoFs at its central vertex, so distinct vertex blocks have no shared unknowns, although
their supporting cells overlap. The overlapping Schwarz patches of Witte et al.~\cite{witte2025tensor} correct a larger
set of unknowns, shared with neighboring patches. Our smoother therefore resolves coupling between vertex blocks
through successive updates, while the reference includes more of this coupling within each local solve. This difference
may contribute to the observed iteration counts.

\subsection{Adaptive Refinement}
\label{sec:biharmonic_results_adaptive}

\begin{figure}[htbp]
    \centering
\begin{subfigure}[b]{0.315\linewidth}
    \centering

    \includegraphics[width=0.9\linewidth]{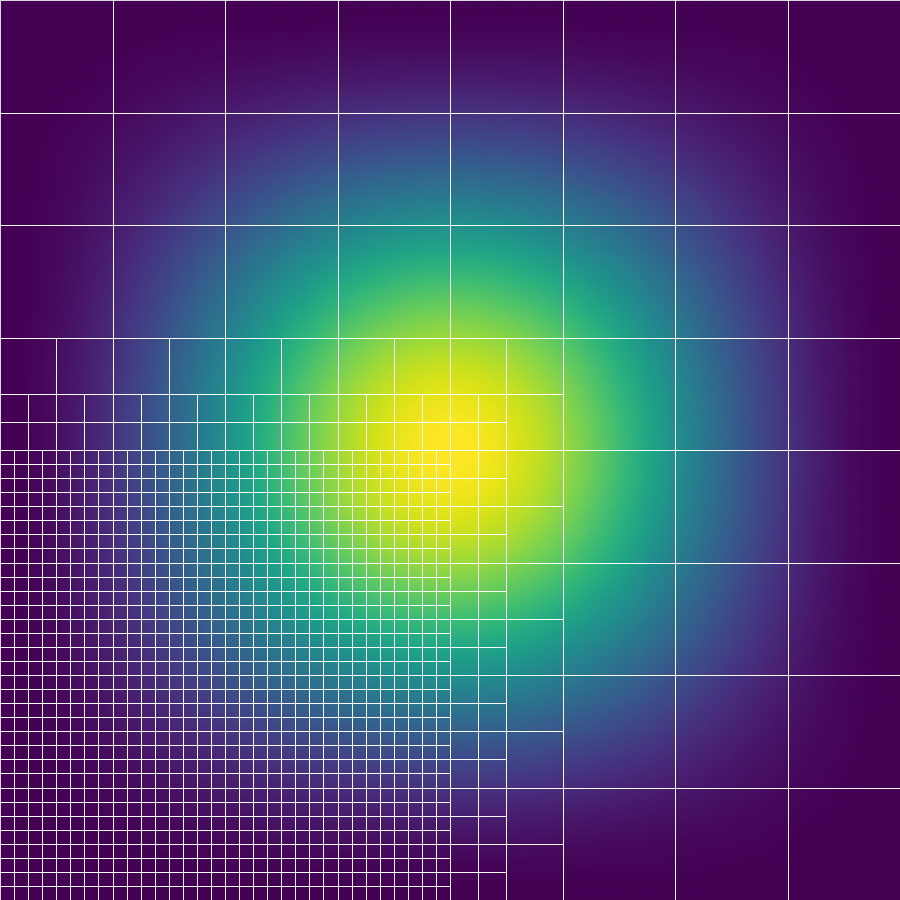}
    \caption{Corner octant.}
    \label{fig:biharmonic_mesh_octant}
\end{subfigure}
\hfill
\begin{subfigure}[b]{0.315\linewidth}
    \centering
    \includegraphics[width=0.9\linewidth]{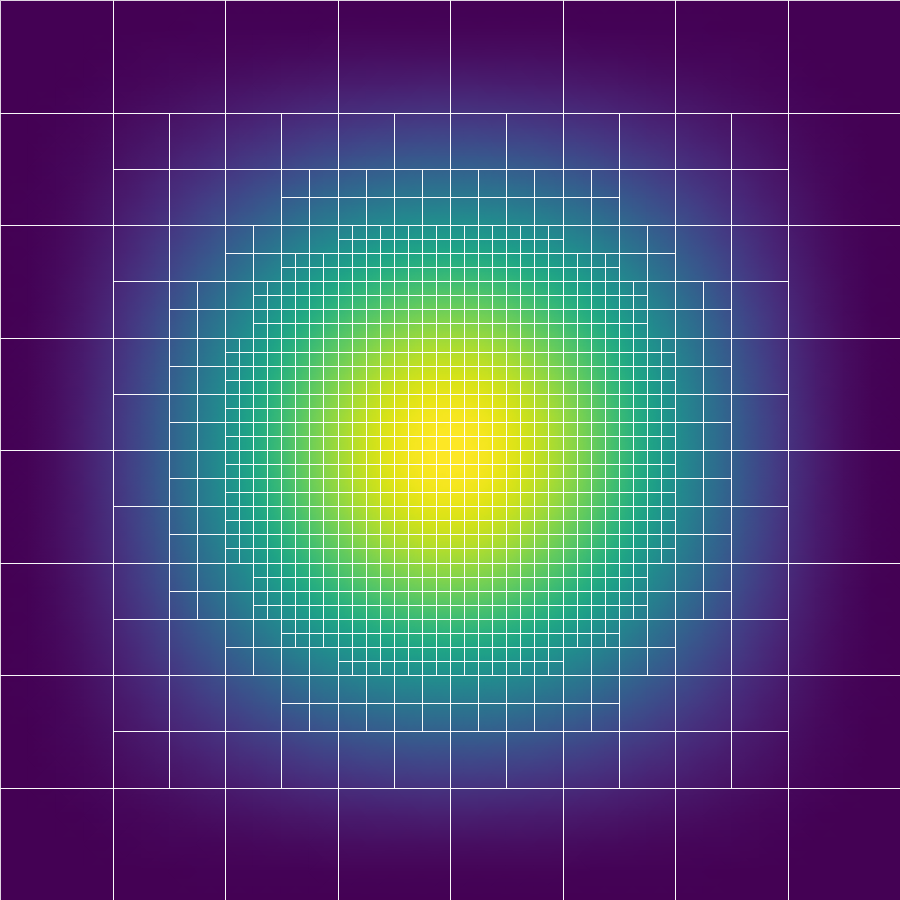}
    \caption{Ball.}
    \label{fig:biharmonic_mesh_ball}
\end{subfigure}
\hfill
\begin{subfigure}[b]{0.315\linewidth}
    \centering
    \includegraphics[width=0.9\linewidth]{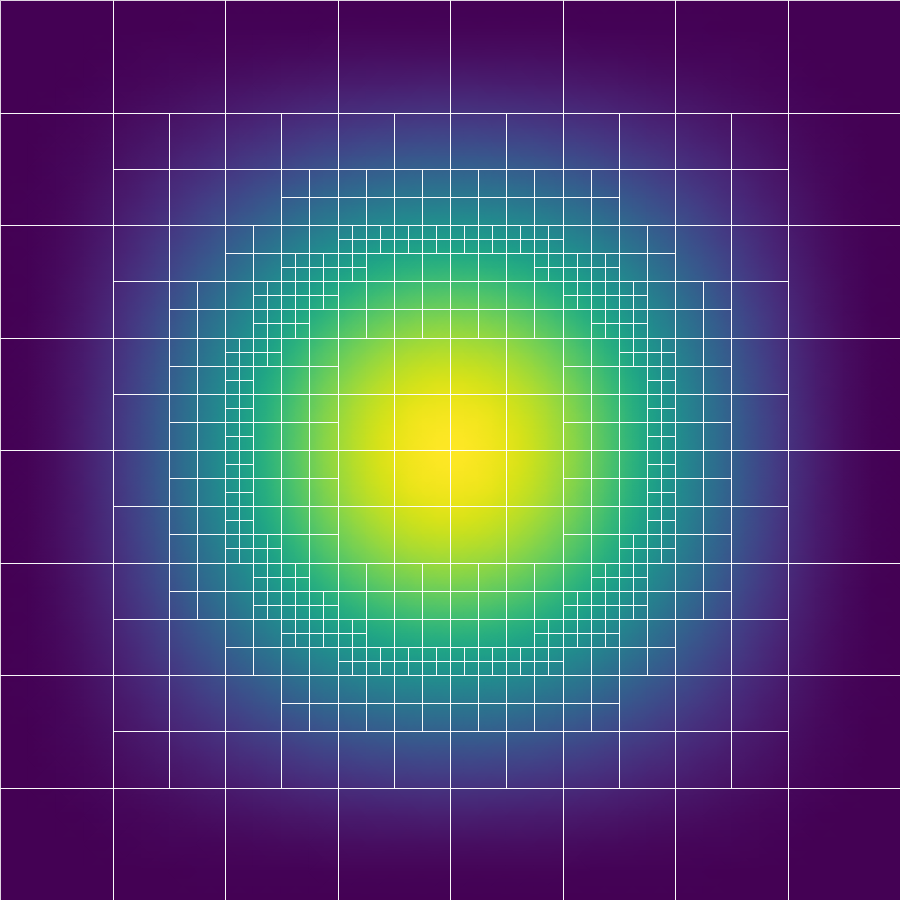}
    \caption{Ring.}
    \label{fig:biharmonic_mesh_ring}
\end{subfigure}
     \caption{Slices at $z=0.45$ through the three locally refined regions at equal adaptive depth. The grid was refined uniformly $n_{\mathrm{uni}}=3$ times, followed by $n_{\mathrm{adap}}=3$ local refinements. The color indicates the computed solution, an approximation of the manufactured solution of \secref{sec:biharmonic_verification} with $L^2$ error $1.007\times 10^{-4}$ for the octant, $1.067\times 10^{-4}$ for the ball, and $1.057\times 10^{-4}$ for the ring.}
    \label{fig:biharmonic_adaptive_meshes}
\end{figure}

Each mesh is identified by the number $n_{\mathrm{uni}}$ of uniform refinements followed by the number
$n_{\mathrm{adap}}$ of nested local refinements; $n_{\mathrm{adap}}=0$ denotes the uniform reference mesh. We use three
marked regions: the corner octant $(0,1/2)^d$, the ball of radius $1/4$ about the domain center, and the cells crossed
by the sphere of that radius. The last region, denoted the ring, forms a shell approximately one and a half cells
thick. Figure~\ref{fig:biharmonic_adaptive_meshes} shows a slice through each region. Each sequence starts from
$n_{\mathrm{uni}}=5$ or $6$ and continues up to the device-memory limit.

All three regions produce refinement edges at which both trace layers are constrained by the masking and transfer
operators of \secref{sec:biharmonic_edge_masking} and \secref{sec:biharmonic_transfer_operators}. Each level consists
of cells of one size, so its cell operator and interior patch inverse can each be reused throughout that level. The
patch admissibility rule excludes updates to dependent DoFs.

Table~\ref{tab:biharmonic_adaptive} reports all three sequences. The errors remain close to, and below, those of the
initial uniform meshes, with limited improvement under further local refinement. The marked regions are fixed in
advance rather than selected by an error indicator, and the error of the smooth manufactured solution is spread over
the whole domain. The squared $L^2$ error is a sum of cell contributions; local refinement reduces only the share of the
marked region, while the unrefined cells keep their mesh size and their contribution. The error therefore levels off at
the contribution of the unrefined cells. These sequences test the solver on refinement edges of controlled geometry,
not the efficiency of adaptive refinement. Continuity across refinement edges holds by construction: the smoothers leave the
value and normal-derivative layers of $E_\ell$ unchanged, and by~\eqref{eq:biharmonic_constraint_prolongation}
prolongation places the interpolated coarse trace data on them. The errors against the manufactured solution, which remain at or
below those of the uniform meshes, check the adaptive discretization end to end.

\begin{table}[tb]
    \centering
    \adjustbox{max width=\textwidth}{%
\small
\begin{tabular}{@{}l cc|rr|rr|rrr@{}}
    \toprule
                     &                       &                       &           &                     & \multicolumn{2}{c|}{Chebyshev--Jacobi} & \multicolumn{3}{c}{Patch, sym.\ + CG}                        \\
    Marked region    & {$n_{\mathrm{uni}}$} & {$n_{\mathrm{adap}}$} & {Cells}   & {$\|u-u_h\|_{L^2}$} & {CG it.} & {[Mcells/s]}                & {CG it.} & {$\rho$} & {[Mcells/s]} \\
    \midrule
    Corner octant & 5 & 0 &    32,768 & $4.924\cdot10^{-7}$ &  6 & 0.580 &  6 & 0.048 & 0.646 \\
                  &   & 1 &    61,440 & $4.598\cdot10^{-7}$ &  7 & 0.803 &  8 & 0.105 & 0.794 \\
                  &   & 2 &   296,535 & $4.451\cdot10^{-7}$ &  7 & 2.45  &  8 & 0.104 & 1.89 \\
                  &   & 3 & 2,153,726 & $4.450\cdot10^{-7}$ &  7 & 3.12  &  8 & 0.105 & 2.50  \\
    \cmidrule(l){2-10}
                  & 6 & 0 &   262,144 & $3.078\cdot10^{-8}$ &  5 & 3.24  &  6 & 0.044 & 2.49  \\
                  &   & 1 &   491,520 & $2.867\cdot10^{-8}$ &  6 & 2.99  &  7 & 0.073 & 2.33 \\
                  &   & 2 & 2,348,711 & $2.830\cdot10^{-8}$ &  6 & 3.56  &  7 & 0.072 & 2.77  \\
    \midrule
    Ball          & 5 & 0 &    32,768 & $4.924\cdot10^{-7}$ &  6 & 0.575 &  6 & 0.048 & 0.690 \\
                  &   & 1 &    52,032 & $4.574\cdot10^{-7}$ & 11 & 0.453 & 11 & 0.217 & 0.531 \\
                  &   & 2 &   195,112 & $4.557\cdot10^{-7}$ & 11 & 1.23  & 11 & 0.228 & 1.30  \\
                  &   & 3 & 1,246,344 & $4.560\cdot10^{-7}$ & 11 & 1.85  & 11 & 0.228 & 1.79  \\
    \cmidrule(l){2-10}
                  & 6 & 0 &   262,144 & $3.078\cdot10^{-8}$ &  5 & 3.22  &  6 & 0.044 & 2.49  \\
                  &   & 1 &   399,120 & $2.871\cdot10^{-8}$ &  9 & 2.07  &  9 & 0.164 & 1.81  \\
                  &   & 2 & 1,448,000 & $2.944\cdot10^{-8}$ &  9 & 2.26  & 10 & 0.183 & 1.96  \\
    \midrule
    Ring          & 5 & 0 &    32,768 & $4.924\cdot10^{-7}$ &  6 & 0.585 &  6 & 0.048 & 0.624 \\
                  &   & 1 &    41,056 & $4.882\cdot10^{-7}$ & 16 & 0.267 & 16 & 0.370 & 0.311 \\
                  &   & 2 &    83,840 & $4.805\cdot10^{-7}$ & 14 & 0.532 & 14 & 0.319 & 0.609 \\
                  &   & 3 &   256,880 & $4.780\cdot10^{-7}$ & 13 & 1.19  & 13 & 0.290 & 1.37  \\
    \cmidrule(l){2-10}
                  & 6 & 0 &   262,144 & $3.078\cdot10^{-8}$ &  5 & 3.30  &  6 & 0.044 & 2.49  \\
                  &   & 1 &   295,632 & $3.068\cdot10^{-8}$ & 13 & 1.42  & 13 & 0.295 & 1.27  \\
                  &   & 2 &   465,256 & $3.048\cdot10^{-8}$ & 12 & 1.60  & 12 & 0.258 & 1.37 \\
    \bottomrule
\end{tabular}
}
    \caption{Conjugate gradient solves on locally refined meshes with Chebyshev--Jacobi and symmetric vertex-patch
        multigrid preconditioners. The mesh has $n_{\mathrm{uni}}$ initial uniform refinements followed by
        $n_{\mathrm{adap}}$ local refinements; $n_{\mathrm{adap}}=0$ is the uniform reference.
        \emph{Cells} counts active cells, and throughput is $N_{\mathrm{cells}}/t_{\mathrm{solve}}$.
        The $L^2$ error is reported for the Chebyshev--Jacobi solve; $\rho$ is the effective contraction factor of
        the patch-preconditioned CG iteration. Both solvers reduce $\langle\tilde r,z\rangle^{1/2}$ by $10^{-8}$.}
    \label{tab:biharmonic_adaptive}
\end{table}

The largest increase in iteration count occurs at the first local refinement. For $n_{\mathrm{uni}}=5$, subsequent
Chebyshev--Jacobi counts remain seven for the octant and eleven for the ball, while the ring count decreases from
sixteen to thirteen. The effective contraction factors of the conjugate-gradient-accelerated patch solver at the
largest tested depth are $0.105$, $0.228$, and $0.290$, respectively, compared with $0.048$ on the initial uniform
mesh. Across both initial resolutions, the first refinement adds one Chebyshev--Jacobi iteration for the octant, four
or five for the ball, and eight or ten for the ring. The two preconditioners yield iteration counts that differ by at
most one on every adaptive mesh and agree on all five adaptive ring meshes. These similar increases motivate the
comparison with refinement-edge geometry in \secref{sec:biharmonic_results_edge_cost}.

Along the octant sequence from $n_{\mathrm{uni}}=5$, the active cell count grows by a factor of $66$ while the
Chebyshev--Jacobi cost per iteration grows by a factor of ten, and solve throughput rises from $0.580$ to
$3.12$\,Mcells/s. Its largest octant mesh reaches $3.56$\,Mcells/s, $16\%$ below the $4.23$\,Mcells/s that
\secref{sec:biharmonic_results_solve} reports for a uniform mesh of comparable size. The adaptive solve requires six
iterations, compared with five for that uniform mesh. The patch solver reaches $2.50$\,Mcells/s for the largest octant
mesh from $n_{\mathrm{uni}}=5$ and $2.77$\,Mcells/s from $n_{\mathrm{uni}}=6$, approximately $20\%$ and $22\%$ below
the corresponding Chebyshev--Jacobi rates. The patch solver requires eight and seven iterations, respectively, compared
with seven and six for Chebyshev--Jacobi. Over the ring sequence from $n_{\mathrm{uni}}=5$, the
conjugate-gradient-accelerated patch solver costs slightly less per iteration than the Chebyshev--Jacobi solver at
every depth. These timings show that iteration count is a substantial contributor to the solve-time differences among
the tested meshes.

\subsection{Refinement-Edge Geometry and Iteration Counts}
\label{sec:biharmonic_results_edge_cost}

The iteration counts increase with the fraction of level cells incident to refinement edges in the three tested
regions. At the first adaptive step from $n_{\mathrm{uni}}=5$, these fractions are approximately $9\%$ for the octant,
$20\%$ for the ball, and $74\%$ for the ring (Table~\ref{tab:biharmonic_level_geometry}). The corresponding increases
in Chebyshev--Jacobi iteration count are one, five, and ten. The ordering of the iteration counts persists at both
initial resolutions and all tested adaptive depths.

\begin{table}[htbp]
    \centering
\begin{tabular}{@{}l rrr@{}}
    \toprule
    Marked region & edge faces per cell & cells on a refinement edge & admissible patch centers \\
    \midrule
    Corner octant & 0.09 &  $9.1\%$ & $90.9$--$98.8\%$ \\
    Ball          & 0.29 & $20.2\%$ & $86.2$--$98.2\%$ \\
    Ring          & 1.00 & $74.5\%$ & $50.0$--$50.1\%$ \\
    \bottomrule
\end{tabular}
     \caption{Geometry of the local levels for the three marked regions. Edge faces per cell and the fraction of cells
        carrying such a face are given at the first adaptive step from $n_{\mathrm{uni}}=5$; both fall as the level
        grows, and the ordering of the three regions holds at every tested depth. The measurements use a reconstruction
        of the marked regions that reproduces the first-step active cell counts exactly for the octant and ring and
        to within $8\%$ for the ball. It omits $2\!:\!1$ balancing at greater depth.}
    \label{tab:biharmonic_level_geometry}
\end{table}
\FloatBarrier

This ordering is consistent with the surface-to-volume ratios of the marked regions. Although the ball and octant have
the same extent along each coordinate axis, their interface areas and volumes differ. The ball has surface area
approximately $0.79$ and volume $0.065$, while the octant has interior interface area $0.75$ and volume $0.125$. The
ball therefore has approximately twice the interface area per unit volume. The ring has a larger ratio because its
inner and outer interfaces enclose a thin shell. These geometric differences set how much of each level carries
edge-constrained data.

\section{Conclusion}
\label{sec:biharmonic_conclusion}

The cell-wise $C^1$ Hermite formulation yields the expected fourth-order $L^2$ convergence for the clamped
biharmonic problem with cubic elements. Measured jumps in the value and normal derivative vanish on matching
interfaces. At refinement edges, the formulation constrains both trace layers and scales reference derivatives
under geometric prolongation.

The two-axis even--odd decomposition gives the highest throughput among the tested operator evaluations,
reaching $85\%$ of peak memory bandwidth. It splits the element matrix into four dense $16\times16$ sector blocks.

On the finest uniform mesh, Chebyshev--Jacobi preconditioned conjugate gradients require five iterations and reach
$4.23$\,Mcells/s. In cells per second, this is $18$ times the double-precision rate reported for the tensor-product
patch multigrid solver of Witte et al.~\cite{witte2025tensor} at the same polynomial degree. Mixed precision is not
explored here; against their mixed-precision rate, the speedup is a more modest factor of $10$. The symmetric patch
preconditioner requires six iterations and reaches $3.16$\,Mcells/s. The patch sweep requests $8.3$ times the field-data volume of one volumetric operator evaluation. Its arithmetic
and memory-access costs cannot be separated by the full-store control experiment alone.

The largest increase in iteration count occurs at the first local refinement. Further refinement produces nearly
constant counts for the octant and ball and decreasing counts for the ring over the tested depths. Both
preconditioners yield counts that differ by at most one on the adaptive meshes. The ordering across the three
regions is consistent with their refinement-interface area relative to volume; the experiments do not isolate this
ratio from other geometric differences. On the largest octant mesh, Chebyshev--Jacobi and patch preconditioning
reach $3.56$ and $2.77$\,Mcells/s, respectively.
 
\ifarxiv
    \section*{Use of generative AI}
\fi
During the preparation of this work the author used Claude (Anthropic) to assist with drafting and editing the text.
The author reviewed and edited all content and takes full responsibility for the publication.

\ifarxiv
    \bibliographystyle{siam}
\fi

\end{document}